\documentclass[reqno]{amsart}
\theoremstyle{plain}
\usepackage{latexsym}
\usepackage{amssymb, xcolor}
\usepackage{graphicx}
\usepackage{enumerate}
\usepackage{array}
\usepackage{lipsum}
\usepackage{epsfig}
\usepackage{enumitem}
\usepackage[sc]{mathpazo}
\def\XXint#1#2#3{{\setbox0=\hbox{$#1{#2#3}{\int}$}
     \vcenter{\hbox{$#2#3$}}\kern-.5\wd0}}

\newcommand{\Ub}{\mathbb{U}}

\newcommand{\Tb}{\mathbb{T}}

\newcommand{\al}{\mathfrak{a}}

\newcommand{\Fk}{\mathfrak{F}}

\newcommand{\lt}{\left}
\newcommand{\rt}{\right}
\newcommand{\nl}{\newline}

\newcommand{\nn}{\nonumber}

\newcommand{\lm}{\lambda}

\newcommand{\qd}{\quad}

\newcommand{\ep}{\epsilon}

\newcommand{\wt}{\widetilde}

\newcommand{\AI}{\mathcal{A}}
\newcommand{\PPI}{\mathcal{P}}

\newcommand{\SI}{\mathcal{S}}

\newcommand{\MI}{\mathcal{M}}
\newcommand{\VI}{\mathcal{V}}

\newcommand{\OI}{\mathcal{O}}

\newcommand{\QI}{\mathcal{Q}}
\newcommand{\LI}{\mathcal{L}}
\newcommand{\KKI}{\mathcal{K}}
\newcommand{\WI}{\mathcal{W}}

\newcommand{\ti}{\tilde}
\newcommand{\la}{\langle}
\newcommand{\ra}{\rangle}

\newcommand{\BI}{\mathcal{B}}

\newcommand{\R}{\mathrm {I\!R}}
\newcommand{\K}{\mathcal{K}}

\newcommand{\ca}[1]{\mathrm{Card}\lt(#1\rt)}
\newcommand{\dia}{\diamondsuit}

\newcommand{\red}[1]{{\textcolor{black}{#1}}} 
\newcommand{\blue}[1]{{\textcolor{black}{#1}}} 
\newcommand{\bblue}[1]{{\textcolor{black}{#1}}} 
\newcommand{\bbblue}[1]{{\textcolor{black}{#1}}} 
\newcommand{\green}[1]{{\textcolor{black}{#1}}} 
\newcommand{\ggreen}[1]{{\textcolor{black}{#1}}} 
\newcommand{\rred}[1]{{\textcolor{black}{#1}}} 
\newcommand{\rrred}[1]{{\textcolor{black}{#1}}}
 \newcommand{\zred}[1]{{\textcolor{black}{#1}}} 
\newcommand{\org}[1]{{\textcolor{black}{#1}}}
\newcommand{\pblue}[1]{{\textcolor{black}{#1}}} 
\newcommand{\rblue}[1]{{\textcolor{black}{#1}}} 
\newcommand{\xred}[1]{{\textcolor{black}{#1}}} 
\newcommand{\xgreen}[1]{{\textcolor{black}{#1}}}

\newtheorem{a1}{Lemma}
\newtheorem{a2}[a1]{Theorem}

\newtheorem{a5}[a1]{Proposition}
\newtheorem{a6}[a1]{Corollary}
\newtheorem{deff}[a1]{Definition}

\theoremstyle{remark}
\newtheorem{remark}{Remark}

\begin{document}
\title[\rred{On the }Rank-$1$ convex hull of a set arising from a hyperbolic system]{\rred{On the }Rank-$1$ convex hull of a set arising from a hyperbolic system of Lagrangian elasticity}
\author[A. Lorent, G. Peng]{Andrew Lorent, Guanying Peng}
\address{A. L.\\Mathematics Department\\University of Cincinnati\\2600 Clifton Ave.\\ Cincinnati, OH 45221; 
G. P.\\Department of Mathematics\\University of Arizona\\617 N. Santa Rita Ave.\\ Tucson, AZ 85721.
}
\email{lorentaw@uc.edu, gypeng@math.arizona.edu.}
\maketitle

\begin{abstract} We address the questions (P1), (P2) asked in \cite{kms} \blue{concerning} the structure of the Rank-$1$ convex hull of \blue{a} submanifold $\KKI_1\subset M^{3\times 2}$ that \blue{is related to} weak solution\blue{s} of the two by two system of \blue{Lagrangian equations of elasticity} studied by DiPerna \cite{dp} with one entropy augmented. This system serves as a model problem for higher order systems for which there are only finitely many entropies. The Rank-$1$ convex hull is of interest in the study of solutions via convex integration: the Rank-$1$ convex hull needs to be sufficiently non-trivial for convex integration to be possible. Such non-triviality is typically \blue{shown} by embedding a $\Tb_4$ (Tartar square) into the set; see for example \cite{mulsv2}, \cite{mulsvri}. We show that in the strictly hyperbolic, genuinely nonlinear case considered by DiPerna \cite{dp},  no $\Tb_4$ configuration can be embedded into $\KKI_1$. 
\end{abstract}

\section{Introduction}

There has recently been a lot of progress on a number of outstanding problems in PDE by reformulating the PDE as a differential inclusion.  In \cite{mulsv1} \rblue{counterexamples} to partial regularity of weak solutions to elliptic systems that arise as the critical point of a strongly quasiconvex 
functional were provided\,\footnote{Contrast this with the well known result of Evans \cite{evans} that minimizers do have partial regularity.}. This was later extended to polyconvex 
functionals in \cite{sz2} and parabolic systems in \cite{mulsvri}. Prior to this Scheffer \cite{sch} provided 
\rblue{counterexamples}  to related regularity problems. In \cite{camles1}, \rblue{De Lellis and Sz\'{e}kelyhidi} reproved \ggreen{(and considerably strengthened)} the well known result of \ggreen{Scheffer \cite{sch2}} on weak solutions to the Euler equation with compact support in space and time, with a much shorter and simpler proof via reformulation as a  differential inclusion. \ggreen{ Previously  Shnirelman \cite{shnir} provided a somewhat simpler proof by a different method}. The advance provided by \cite{camles1} opened an approach to Onsager's conjecture which was subsequently studied intensively by a number of authors \cite{camles3},  \cite{camles2}, \cite{buvi2}, \cite{phil}, \cite{phil2} with a final solution being provided by \cite{phil1}, \cite{camles8}. \rrred{Further work brings these methods} to the study of the Navier-Stokes \rred{equations} \cite{buvi}. An excellent recent survey is provided by 
\cite{camlas10}. The general term used to describe the method of constructing solutions of PDE via differential inclusions is \em convex integration\rm. Indeed the antecedent to many of these results are the celebrated results of Nash \cite{nash}, Kuiper \cite{kuip} and Gromov \cite{grom}. \xred{Recently, De Lellis,  De Philippis, Kirchheim and Tione \cite{dede} studied the question of regularity of 
	\em stationary points \rm of strictly polyconvex functionals from the perspective of (generalized) differential inclusions and convex integration, as will be explained - their results (while more general and for a completely different problem and differential inclusion)  are in spirit quite similar to ours.}

The purpose of this paper is to contribute to the study of regularity and uniqueness of \ggreen{\em entropy solutions\rm} of systems of conservation laws  via differential inclusions and convex integration. \ggreen{By this we mean solutions that satisfy (in a distributional sense)
	entropy inequalities of the form $(\eta(u))_t+(q(u))_x\leq 0$ for \bbblue{all} entropy/entropy-flux pair\bbblue{s} $(\eta,q)$; see Definition (36), (37) in Section 11.4, \bbblue{\cite{evans3}}}. The first step in such a program is to consider a PDE \ggreen{and adjoined entropy inequalities} reformulated  as a differential inclusion into a submanifold $\KKI\subset M^{m\times n}$ (\blue{the set of $m\times n$ matrices}) and to determine if $\KKI$ admits a four matrix configuration known as $\Tb_4$ configuration, or \em Tartar square\rm\,\footnote{Indeed as noted in \cite{mulsv2}, $\Tb_4$ configurations played an important role in \cite{sch} and seem to have been discovered independently by a number of authors.}. We will describe this configuration and its $n$-matrix variants in more detail in Section \ref{S2.2}. We study a simple two by two system that arises from the Lagrangian formulation of elasticity and is augmented \blue{by} one entropy/entropy flux pair. This system can be reformulated as a differential inclusion into a submanifold $\KKI_1\subset M^{3\times 2}$. The study of this system and its associated submanifold $\KKI_1$ was initiated by Kirchheim, M\"{u}ller, \v{S}ver\'{a}k in 
\cite{kms}, Section 7. They provided a hierarchy of \blue{properties} (P1), (P2), (P3), (P4) and asked \blue{for the hypotheses on the system under which (P1)--(P4) hold}. In \cite{lope} we investigated the system and answered the question on (P4). Non-technically speaking, \blue{the properties} (P1)--(P4) \blue{concern} a hierarchy of hulls of $\KKI_1$. Non-triviality of the hull associated with (P1), \blue{(P2)} (\blue{the \em Rank-$1$ convex hull \rm of $\KKI_1$}) would open the 
prospect of an infinity of solutions to the differential inclusion into $\KKI_1$. The hull associated with \blue{(P3)}, (P4) (the \em polyconvex hull \rm of $\KKI_1$) contains the Rank-$1$ convex hull of $\KKI_1$ and the result of \cite{lope} (see Section \ref{S2.2}) - specifically \rred{that} the polyconvex hull is non-trivial \rred{when} the system \rred{is hyperbolic} -  opened the possibility 
that the structure of $\KKI_1$ \blue{is} sufficiently rich to allow for an infinity of solutions to the differential inclusion into $\KKI_1$. The Rank-$1$ convex hull would be non-trivial if a $\Tb_4$ configuration could be found in $\KKI_1$. \blue{Unfortunately} we show in this paper that no $\Tb_4$ exists in $\KKI_1$ \rred{when the system is hyperbolic and genuinely nonlinear in the sense of DiPerna \cite{dp}} (\blue{see} Theorem \ref{T1}). This does not rule out the possibility of embedding $n$-matrix version of $\Tb_4$ (denoted by $\Tb_n$) in $\KKI_1$ (as for example was shown in \cite{sz2} for $\Tb_5$) and non-triviality of the Rank-$1$ convex hull of $\KKI_1$. However, in establishing non-triviality of the Rank-$1$ convex hull of a set, \bblue{an} \blue{important} first step is to \bblue{understand the possibility of embedding $\Tb_4$ configurations inside the set; see \cite{kms}, Section 3.5, where non-existence of $\Tb_4$ configurations in an important setting is proved}, and \rblue{\cite{ta3}, Remark 10 (also \cite{kms}, Proposition 19), \cite{kms}, Proposition 21 and \cite{sz}} for close connections between non-triviality of the Rank-$1$ convex hull and existence of $\Tb_4$ configurations in certain sets without Rank-$1$ connections.  \blue{For this reason} we complete this study of $\Tb_4$ configurations for the set $\KKI_1$.  \xred{In a recent work \cite{dede}, the authors formulated a more general kind of differential inclusion (that they named a \em div-curl \rm differential inclusion) into set of matrices $K_f$ \rblue{whose solution corresponds to a Lipschitz \em stationary point \rm \rblue{of the energy $\int_{\Omega}f(Du)\,dx$  for polyconvex functions $f$}. Their main result (\cite{dede}, Theorem 1.2) established the strong result that if $f\in C^1(\R^{n\times m})$ is strictly polyconvex then $K_f$ does not contain a $\rblue{\Tb'_N}$ configuration for any $N$, where  $\rblue{\Tb'_N}$ is a generalization of $\Tb_N$ adapted to \em div-curl \rm differential inclusions. This result is an important step towards the long open problem of regularity of stationary points of \rblue{the above energy} via convex integration solutions of the \em div-curl \rm  differential inclusion into $K_f$.} }

	\subsection{Conservation laws}
	
	A scalar conversation law \blue{in space dimension one for an unknown function $u(x,t)$} is an equation of the form 
	\begin{equation}
		\label{eqinta1}
		u_t+\lt(f(u)\rt)_x=0.
	\end{equation}
	It is not hard to see there are infinitely many weak solutions. To select the \blue{physically} correct 
	solution, the notion of \em entropy/entropy flux \rm pair was introduced. This is a pair of functions $(\eta,q)$ where $\eta$ is convex and $q'=\eta' f'$. If $u$ is a 
	smooth solution to (\ref{eqinta1}) we have that $\lt(\eta(u)\rt)_t+\lt(\blue{q}(u)\rt)_x=0$. If we regularize the equation (\ref{eqinta1}) by forming $u^{\ep}_t+\lt(f(u^{\ep})\rt)_x=\ep u^{\ep}_{xx}$,  then assuming 
	$\{u^{\ep}\}_{\ep>0}$ is bounded in $L^{\infty}\lt(\R\times \lt(0,\infty\rt)\rt)$, the method of compensated compactness (see \cite{evans2}, Chapter 5, Section D) allows us to conclude that 
	$u^{\ep}\overset{L^1}{\rightarrow} u$ for some weak solution $u$ of (\ref{eqinta1}). Further it turns out that $\mathrm{div}(\eta(u),q(u)):=\lt(\eta(u)\rt)_t+\lt(q(u)\rt)_x$ forms a negative measure for every entropy/entropy flux pair $(\eta,q)$. We call solutions of (\ref{eqinta1}) that satisfy this property \em entropy solutions\rm. For scalar conservation laws \blue{at least in \rblue{one space dimension}} this is the correct notion, namely, entropy solutions enjoy uniqueness, regularity and can even be described in closed form for sufficiently regular $f$; see \cite{evans3}, Theorem 3 in Section 11.4 and \cite{ol}, Section 3.4.2.
	
	The theory for systems of conservation laws \ggreen{in one space dimension} is much more limited. The two \ggreen{main} methods to produce existence of solutions are Bressan's semigroup method for (small) BV initial data \cite{bres}, \cite{bres2} and the compensated compactness method pioneered by Tartar, Murat and DiPerna \cite{ta1}, \cite{ta2}, \cite{murat}, \cite{dp2}, \cite{dp} and developed by many others. The compensated compactness method proceeds by finding appropriate entropies for the system under consideration and under reasonable assumptions on a regularizing sequence, proving compactness and hence existence of $L^{\infty}$ solutions that satisfy an entropy \blue{production} inequality of an analogous form to the scalar equation. \ggreen{Indeed if we expect the ``physically correct" solution to a system of conservation laws to be the limit of solutions $u^{\ep}$ to the system with an additional viscosity term $\ep u^{\ep}_{xx}$, assuming compactness can be established as $\ep\rightarrow 0$, then the limiting function $u$ will be an entropy solution; see Theorem 2 in Section 11.4, \bbblue{\cite{evans3}}. For this reason and the fact that it is the correct notion for scalar conservation laws, we \org{are interested to }study the question of uniqueness and regularity of 
		entropy solutions of systems of conservation laws in one space dimension.} 
	
	Given the success of the method of convex integration in addressing \bbblue{related} questions for elliptic systems, the Euler equation and the Navier-Stokes equation, a natural goal (already implicit in \cite{kms}) is to extend the scope of such approach to construct \rblue{counterexamples} to \ggreen{uniqueness and regularity} for systems of conservation laws\,\footnote{This goal and this approach has been introduced to us by V. \v{S}ver\'{a}k \cite{svper}.}. 
	
	The system chosen for study in \cite{kms} is the two by two system of Lagrangian equations of elasticity given by
	\begin{equation}
		\label{eqint50}
		\begin{cases}
			v_t - u_x=0,\\
			u_t - \al(v)_x = 0
		\end{cases}
	\end{equation}
	\blue{for the unknowns $u, v$ and some appropriate function $\al$. This system was studied earlier by DiPerna \cite{dp2}, \cite{dp} under the assumption that $\al'>0$, i.e., the system is hyperbolic and additional assumptions on the sign of $\al''$. In \cite{dp2}, DiPerna proved existence of solutions to the system (\ref{eqint50}) using the method of compensated compactness with the help of all entropy/entropy flux pairs. Possibly motivated by the question of compactness for higher dimensional systems, in \cite{dp}, he proved a local existence result when the system is genuinely nonlinear, i.e., $\al''\ne 0$ with just two physical entropy/entropy flux pairs}. Following \cite{dp} we introduce the natural entropy/entropy flux pair $(\eta_1, q_1)$ defined by
	\begin{equation*}
		\eta_1(u,v) := \frac{1}{2}u^2+\Fk(v), \qd q_1(u,v) := -u\al(v),
	\end{equation*}
	\blue{where $\Fk$ is an antiderivative of the function $\al$. Another dual entropy/entropy flux pair $(\eta_2,q_2)$ was also introduced in \cite{dp}. We omit the technical formulas for the dual pair since it is not relevant in this paper. The results in \cite{dp} \rred{demonstrate} that the system (\ref{eqint50}) augmented by the two entropy/entropy flux pairs $(\eta_i,q_i)$ is rigid enough for the method of compensated compactness to work. A natural question is to further understand this system coupled with just one entropy/entropy flux pair, and in particular, to understand the uniqueness of solutions. For higher order systems, there are only finitely many entropy/entropy flux pairs, and thus it is of great importance to understand the structure of systems augmented by only a few entropy/entropy flux pairs. For this reason, the system (\ref{eqint50}) coupled with $(\eta_1,q_1)$ serves as a model problem and was singled out in \cite{kms}.} 
	
	\blue{As in \cite{kms}, we consider weak solutions $(u,v)$ of the following system
		\begin{equation}\label{eq203}
			\begin{cases}
				v_t - u_x=0,\\
				u_t - \al(v)_x = 0,\\
				(\eta_1(u,v))_t + (q_1(u,v))_x \leq 0.
			\end{cases}
		\end{equation}
		This system can be formulated as a differential inclusion into the set\,\rblue{\footnote{Note that a differential inclusion into set $\K_1$ gives a solution to (\ref{eq203}) with the inequality replaced by an equality.}} $\K_1$ given by
		\begin{equation}\label{ep100}
			\KKI_1 := \lt\{\lt(\begin{matrix}
				u & v\\
				\al(v) & u\\
				u\al(v) & \frac{1}{2}u^2+\Fk(v)\\
			\end{matrix}\rt): u,v\in\R\rt\}.
		\end{equation}
		(See \cite{kms}, Section 7 for the details.) For the convenience of later discussions, we define $P:\R^2\rightarrow M^{3\times 2}$ by
		\begin{equation}
			\label{eq4}
			P(u,v):=\lt(\begin{matrix}
				u & v\\
				\al(v) & u\\
				u \al(v) & \frac{1}{2}u^2+\Fk(v)\\
			\end{matrix}\rt). 
		\end{equation}
		If there is a way to construct convex integration solutions to the differential inclusion into the set $\K_1$, \blue{a consequence would be} non-uniqueness of solutions to (\ref{eq203}). The construction of the former would require the Rank-$1$ convex hull of $\K_1$ to be sufficiently large. For this reason, the questions raised in \cite{kms} concern the various hulls of the set $\K_1$ and we will discuss this in more detail in the next subsection.} 
	
	\subsection{Convex integration, Tartar squares, Rank-$1$ convex and Polyconvex hulls}
	\label{S2.2} A basic building block for \blue{non-trivial solutions to a} differential inclusion is the existence of \em Rank-$1$ connections \em within a set $\K$. We say $A,B\in \K$ are \em Rank-$1$ connected \rm if $\mathrm{Rank}(A-B)=1$. Restricting to $\K\subset M^{2\times 2}$ for simplicity\,\footnote{For the general case in $M^{m\times n}$ the construction is the same, simply slightly harder to visualize.},  we see that 
	$A,B$ are Rank-$1$ connected if and only if there exists some $v\in S^1$ such that $Av=Bv$. By cutting a square with sides parallel to $v$ and $v^{\perp}$ into strips parallel to $v$, \blue{we can construct a Lipschitz mapping $u$ with $Du$ taking the values $A$ and $B$ alternately in adjacent strips. This mapping $u$ satisfies the differential inclusion $Du\in\{A,B\}$} and is not affine, and is referred to as a \em laminate\rm; see \cite{mul}, Section 2.1. Given that this is the most natural way to build a differential inclusion, a natural conjecture might be that if a set $\K$ contains no Rank-$1$ connections then no non-trivial differential inclusion into it can be built. This is false and the 
	first hint as to why comes from the Tartar square or $\Tb_4$ configuration. Identifying diagonal matrices with points in the plane via $\rred{\Pi}:\lt(\begin{smallmatrix} a & 0 \\ 0 & b \end{smallmatrix}\rt)\mapsto 
	\lt(\begin{smallmatrix} a \\  b \end{smallmatrix}\rt)$ we see that diagonal matrices $D_1, D_2$ are Rank-$1$ connected if and only if $\rred{\Pi}(D_1)$ and $\rred{\Pi}(D_2)$ lie on the same vertical or horizontal line. With this in mind it is not hard to see that the set $\K:=\lt\{A_1,A_2,A_3,A_4\rt\}$ given by 
	\begin{equation}
		\label{eqint6}
		A_1=-A_3=\mathrm{diag}(-1,-3)\text{ and }A_2=-A_4=\mathrm{diag}(-3,1)
	\end{equation}
	does not have Rank-$1$ connections. Nevertheless we can construct a sequence $\{u_k\}$ with the property that $\mathrm{dist}(Du_k,\K)\rightarrow 0$ in measure and 
	$Du_k$ does not converge in measure; see Lemma 2.6 in \cite{mul}. 
	
	It turns out that the heart of this is the fact that \blue{the set $\K$ defined above forms a $\Tb_4$ configuration and} the \em Rank-$1$ convex hull \rm of $\K$ is non-trivial. More generally, we give
	\begin{deff}
		\label{def1}
		An ordered set of $N\geq 4$ matrices $\lt\{T_i\rt\}_{i=1}^N \subset M^{m\times n}$ without Rank-$1$ connections 
		is said to form a $\Tb_N$ configuration if there exist matrices $P_i,C_i\in M^{m\times n}$ and numbers $\kappa_i>1$ such that 
		\begin{equation}
			\label{eq1}
			\begin{split}
				T_1&=P+\kappa_1 C_1,\\
				T_2&=P+C_1+\kappa_2 C_2,\\
				&\dots \\
				T_N&=P+C_1+C_2+\dots C_{N-1}+\kappa_N C_N,
			\end{split}
		\end{equation}
		where $\mathrm{Rank}(C_i)=1$ for all $i$ and 
		\begin{equation}
			\label{eq3.5}
			\sum_{i=1}^N C_i=0.
		\end{equation}
	\end{deff}
	\begin{remark}\label{r:nodegenerateT4}
		It is a simple consequence of the above definition that $\Tb_N$ configurations cannot be contained in an affine space of dimension one.
	\end{remark}
	We say a function $f:M^{m\times n}\rightarrow \R$ is Rank-$1$ convex if $f\lt(\lm A+(1-\lm)B \rt)\blue{\leq}\lm f(A)+(1-\lm)f(B)$ whenever $\mathrm{Rank}(A-B)=1$. The  Rank-$1$ convex hull of a \blue{compact} set $\K$ is defined as (see \blue{\cite{kms}, Section 2}) 
	\begin{equation}
		\label{eqint9}
		\K^{rc}:=\lt\{F\in M^{m\times n}: f(F)\leq \blue{\sup_{\K}} f\text{ for all Rank-1  convex } f:M^{m\times n}\rightarrow \R\rt\}. 
	\end{equation}
	\blue{For a general set $E$ we set}
	\begin{equation*}
		\blue{E^{rc}=\bigcup_{\K\subset E \text{ compact}}\K^{rc}}.
	\end{equation*}
	Now \blue{a celebrated result of M\"{u}ller and \v{S}ver\'{a}k} (see Theorem 1.1 in \cite{mulsv3}) states that if $\Omega$ is a Lipschitz domain and $\K\subset M^{m\times n}$ is open \blue{and bounded}, then there exists a solution to the differential inclusion $Du\in \K$ a.e. with $u=v$ on $\partial \Omega$, where 
	$v$ is a \blue{piecewise} affine map with\,\rblue{\footnote{Here we are stating a more restrictive version of their theorem to avoid some technicalities.}} $Dv\in \K^{rc}\backslash \K$. Hence a non-trivial solution to the differential inclusion into $\K$ exists. However for applications to PDE, it is not generally the case that the set $\K$ is open. The proofs of \cite{mulsv2}, \cite{mulsvri} work by showing that many $\Tb_4$ configurations can be embedded into $\K$, specifically 
	$\Tb_4$ configurations that can be perturbed so that the embedded $\Tb_4$ moves in a ``transversal" way. \blue{Although a necessary} condition \blue{for the existence of \rred{(periodic)} non-trivial solutions to a differential inclusion into a set $\K$ is the non-triviality of $\K^{rc}$, the latter is not sufficient} (for example 
	it is known \cite{chki} that there is no \blue{non-trivial} differential inclusion into $\lt\{A_1,A_2,A_3,A_4\rt\}$, where $A_i$ are defined in (\ref{eqint6}), however $\lt\{A_1,A_2,A_3,A_4\rt\}^{rc}\not=\lt\{A_1,A_2,A_3,A_4\rt\}$). Despite this, in many or even most circumstances \blue{non-triviality of $\K^{rc}$} is enough; see for example the recent interesting work on $\Tb_5$ configurations \cite{fosz}. 
	
	Thus with a view to constructing non-trivial differential inclusions into $\KKI_1$ defined in (\ref{ep100}), in \cite{kms} the authors \blue{asked about the condition on the function $\al$ such that $\K_1^{rc}$ is trivial or non-trivial at least locally and this is basically the content of (P1).} With respect to non-triviality this is the hardest of a hierarchy of questions (P1)--(P4). To explain this further we need to introduce some more concepts. Let $\PPI(\K)$ denote 
	the set of probability measures on $M^{m\times n}$ that are supported on $\K$, and given $\nu\in\PPI(\K)$, let $\la \nu,f \ra:=\int f(X) d\nu(X)$ and $\bar\nu$ be the barycenter of $\nu$. Following 
	\cite{kms},  Section 4.2 we define 
	\begin{equation}
		\label{eqint101}
		\MI^{rc}(\K):=\lt\{\mu\in \PPI(\K): \la \mu, f\ra \geq f\lt(\bar\mu\rt) \text{ for all Rank-1 convex functions }f \rt\}. 
	\end{equation}
	One of the most useful characterizations of $\K^{rc}$ \blue{for compact $\K$} is that $\K^{rc}=\lt\{\bar\mu: \mu\in \MI^{rc}(\K)\rt\}$,  see \blue{\cite{kms}, Section 4.2}. A particular very useful subclass of Rank-$1$ convex functions is the set of \em Polyconvex functions\rm, \blue{which can be expressed as convex functions of minors}. The analog to $\K^{rc}$ and $\MI^{rc}(\K)$ (recall (\ref{eqint9}), (\ref{eqint101})) are the polyconvex hull $\K^{pc}$ and the set of probability measures $\MI^{pc}(\K)$ that are defined in exactly the same way but with respect to polyconvex functions. Since polyconvex functions form a strict subclass of Rank-$1$ convex functions, we have the inclusions 
	\begin{equation}
		\label{introeqa1}
		\ggreen{\K^{rc}\subset \K^{pc}\text{ and }\MI^{rc}(\K)\subset \MI^{pc}(\K).} 
	\end{equation}
	In \cite{lope} we named the measures in $\MI^{pc}(\K)$ \em Null Lagrangian measures \rm and studied \zred{necessary and sufficient conditions on subspaces in $M^{m\times n}$ to support non-trivial Null Lagrangian measures and also question (P4) of \cite{kms}. With respect to the latter, we showed that given $(u_0,v_0)\in\R^2$,} if $\al'(v_0)>0$ (\blue{the system is hyperbolic}) then in any neighborhood $U$ of $P(u_0,v_0)$  (recalling (\ref{eq4})), $\MI^{pc}(U\cap \KKI_1)$ is non-trivial. On the other hand, if $\al'(v_0)<0$ (\blue{the system is elliptic}) then $\MI^{pc}(U\cap \K_1)$ is trivial (the latter case is to be expected). This result opens up the hope that for $\al'(v_0)>0$,  the set $\lt(U\cap \K_1\rt)^{rc}$ could also be non-trivial and 
	a non-trivial differential inclusion into $\K_1$ could be obtained. This would be an important first result in the study of non-uniqueness of \ggreen{entropy solutions} to systems of hyperbolic conservation laws via convex integration. The credit for this question and this formulation belongs to the authors of \cite{kms}. 
	
	\org{Note that the vast  majority of theorems that establish existence of solutions via compensated compactness essentially comes down to showing $\MI^{pc}(\K)$ consists of Dirac measures (assuming appropriate bounds on the \pblue{approximating sequence}) where 
		$\K\subset M^{m\times n}$ is the submanifold defined by the systems and the augmented entropies (just as $\K_1$ is defined by (\ref{eq203})). The only example of compensated compactness that we are aware of that does not proceed by establishing triviality of Null Lagrangian measures is \v{S}ver\'{a}k's proof of compactness for the three well problem \pblue{based on triviality of the \em Quasiconvex hull \rm $\mathcal \K^{qc}$ (see \cite{mul}, Section 4.4; this is sandwiched between $\mathcal \K^{rc}$ and $\mathcal \K^{pc}$)}; see page 298 in \cite{sv5} and Theorem 2.5 in \cite{mul}\,\footnote{It is likely that the sharp results of \cite{faclas} could also be used to generate explicit examples in $M^{2\times 2}$.}. As such for systems for which existence has been established via compensated compactness, (\ref{introeqa1})  implies that the Rank-$1$ convex hull of the set $\K$ is trivial and there is no hope to prove non-uniqueness via differential inclusions and convex integration.}
	
	\ggreen{So given a system of conservation laws \bbblue{augmented by finitely many entropies}, from the perspective of differential inclusions there are essentially two ``levels" at which entropy solutions could be shown to be not a viable notion of solution\,\footnote{\bbblue{The two by two system (\ref{eqint50}) has infinitely many entropies, and it is known from \cite{dp} that the method of compensated compactness works even for the system adjoined by two appropriate entropies.} It seems to the authors of this paper that for two by two systems \bbblue{augmented by} infinitely many entropies there is little hope to \rblue{counterexamples} of uniqueness and regularity by differential inclusions and \bbblue{convex integration}.}. The first and lower level is to show that the set $\K$ (of the associated differential inclusion) supports non-trivial Null Lagrangian measures (i.e. $\mathcal M^{pc}(\K)$ \bbblue{contains measures that are not Diracs}). This \org{means that a proof of triviality of the Quasiconvex hull $\mathcal \K^{qc}$ is required to }construct solutions via compensated compactness methods. \org{Quasiconvex functions are not well understood. Despite some powerful recent advances in $M^{2\times 2}$ \cite{faclas}, from the perspective of conservation laws this would seem to be a very hard (though not impossible) task}. If this first level is reached, a second deeper level is to show that $\K^{rc}$ is sufficiently non-trivial that \bbblue{non-trivial} solutions to the differential inclusion $Dw\in \K$ can be constructed via convex integration. This second level shows that entropy solutions are not the correct notion since in this case solutions are wildly non-unique and have no regularity beyond Lipschitzness. \org{Further if $\K^{rc}$ could merely be shown to be non-trivial, this alone \pblue{wipes out the possibility of establishing} the existence of solutions via compensated compactness since $\K^{rc}\subset \K^{qc}$; see equation (4.8) and Theorem 4.7 in \cite{mul}}. The first level is represented by questions (P3), (P4) of \cite{kms} and questions (P1), (P2) are directed towards the second level.} 
	
	In this paper we make the first progress in answering the questions in (P1), (P2) of \cite{kms} regarding the structure of $\K_1^{rc}$ \blue{by investigating the possibility of embedding $\Tb_4$ configurations in $\K_1$}. If this could be done, an immediate consequence would be the non-triviality of $\K_1^{rc}$. Unfortunately our main result shows that no $\Tb_4$ can be embedded into $\K_1$ under \zred{the assumptions of hyperbolicity and genuine non-linearity (in the sense of DiPerna \cite{dp}) of} the system (\ref{eqint50}). Specifically, we prove
	\begin{a2} 
		\label{T1}
		Suppose \blue{$\al\in C^2(\R)$} is \blue{strictly} increasing and strictly convex, and let the set $\K_1$ be defined in (\ref{ep100}). Then 
		$\KKI_1$ does not contain $\Tb_4$ configurations. 
	\end{a2}
	\begin{remark}
		\blue{With only very minor modifications, our proof of Theorem \ref{T1} also rules out $\Tb_4$ configurations in the set $\KKI_1$ if the function $\al$ is strictly increasing and strictly concave.}
	\end{remark}
	
	\rred{Theorem \ref{T1} easily implies a local version:}
	\begin{a6} 
		\rred{Suppose $\al\in C^2(\R)$ with $\al'(v_0)>0$ and $\al''(v_0)>0$ for some $v_0\in \R$, then for any $u_0$ there exists some neighborhood $U\subset M^{3\times 2}$ of $P(u_0,v_0)$ (defined by (\ref{eq4})) such that $\KKI_1\cap U$ does not 
			contain $\Tb_4$ configurations.}
	\end{a6}
	
	Note that the strict sign condition on $\al''$ is a sufficient condition to rule out Rank-$1$ connections in the set $\K_1$; see \zred{Proposition \ref{TCXE} below and for a \em local \rm result for a more general system see Theorem 4.1 in \cite{dp}}. Thus it is also an important condition from the differential inclusion point of view. Note that if $\al''$ changes sign, then generically the set $\K_1$ contains Rank-$1$ connections. Specifically, in Section \ref{ce} we show
	\begin{a5}
		\zred{\label{TCXE} 
			Let $I\subset \R$ be an \bbblue{open} interval and let $\al\in C^2(I)$ satisfy $\al'>0$ on $I$. Let  $P(u,v)$ be defined by (\ref{eq4}) and define 
			\begin{equation}
				\label{eq4039}
				\K^I_1:=\lt\{P(u,v):v\in I, u\in \R\rt\}.
			\end{equation}
			If the function $\al$ has an isolated inflection point in $I$, then $\K^I_1$ contains Rank-$1$ connections. Conversely if $\al$ is either strictly convex or strictly concave on $I$,  then 
			$\K^I_1$ has no Rank-$1$ connections. }
	\end{a5}
	
	\begin{remark}
		\zred{At the end of \cite{dp}, Section 5, DiPerna conjectures that ``the wave cone associated with a system of conservation laws that is not genuinely nonlinear cannot be separated from the constitutive manifold through the introduction of any finite number of entropy forms". For the system \bbblue{(\ref{eqint50}) adjoined by two entropy forms, he remarks in Section 4, Remark 1 and the end of Section 5 that, if $\al$ has one inflection point, then this fact can be easily verified using the calculations of Section 10.} Proposition \ref{TCXE} and its proof can be thought of as a detailed ``exposition/clarification" of these remarks \bbblue{for the system (\ref{eq203})}. Note further that if $\K^I_1$ contains a Rank-$1$ connection, then the laminate construction sketched at the start of 
			Section \ref{S2.2} \bbblue{gives \rblue{counterexample} to uniqueness of the} system (\ref{eq203}).}
	\end{remark}
	
	\begin{remark}
		\bblue{As a consequence of Proposition \ref{TCXE}, if $\al$ is a strictly increasing real analytic function, then the set $\K_1$ associated to the function $\al$ contains Rank-$1$ connections if and only if $\al$ has an inflection point. It is not clear to the authors whether such equivalence holds true for less regular functions $\al$.}
	\end{remark}
	The conclusion in Theorem \ref{T1} is a negative result in that the more exciting direction would be to establish the existence of $\Tb_4$ inside $\K_1$ \blue{under the assumptions that the system (\ref{eqint50}) is hyperbolic and genuinely nonlinear}. \rblue{However our result does not rule out the possibility of $\Tb_N$ configurations inside $\K_1$}. \xred{A well known example  of a \rblue{set} that does not admit an embedded $\Tb_4$ but does have $\Tb_5$ configurations (\rblue{leading to convex integration solutions of the differential inclusion into the set and answering the important question of regularity of \emph{critical points} of polyconvex functionals}) is given \rblue{in} \cite{sz2}. On the other hand, as mentioned previously, in \cite{dede} the authors established the non-existence of $\Tb_N'$ configurations \rblue{in the set $K_f$} for any $N$. \rblue{These \xgreen{things} suggest that} from Theorem \ref{T1} little can be guessed about the existence of $\Tb_N$ configurations in $\K_1$ (under the assumptions of Theorem \ref{T1})\,\footnote{\rblue{On a somewhat related well known result,} it is known that the differential inclusion into \rblue{any} finite set of four matrices \rblue{without} Rank-$1$ connections has no convex integration solutions  \cite{chki}, however \rblue{there exists a set of five matrices without Rank-$1$ connections that admits  convex integration solutions of the corresponding differential inclusion}; see \cite{kir}, Chapter 4, Section 3.}. Nevertheless we believe our methods will aid in the study of this \rblue{question}. }

	\subsection*{Acknowledgments} \blue{The} first author would like \blue{to} thank V. \v{S}ver\'{a}k for many very helpful discussions during a visit to Minnesota in summer of 2018. The idea to study entropy solutions of system\blue{s} of conservation laws via differential inclusions and convex integration is from him. Also \blue{a} number of key ideas used in \blue{this} paper (in particular Lemmas \ref{L13} and \ref{L12}) are from \v{S}ver\'{a}k \cite{svper2}.  The first author also gratefully acknowledges the support of the Simons foundation, collaboration grant  \#426900. \xred{Both authors warmly thank the anonymous referee for very careful reading of the paper and for pointing out a number of improvements.} Finally we would like to thank Sam Krupa and L\'{a}szl\'{o} Sz\'{e}kelyhidi for pointing out the original proof of the 2d case contained a calculation error, this preprint version contains the rearrangement of the order of lemmas and different arguments required to fix this error.

\section{Sketch of proof}
\label{proofsk}


Let $\K:=\lt\{T_0,T_1,T_2,T_3\rt\}\subset\K_1$ (\blue{this labeling is more convenient for the proofs}) where $T_i=P(u_i, v_i)$ and the mapping $P$ is given  in (\ref{eq4}). Denoting $V_k=T_k-T_0$ for $k=1,2,3$, our first observation is
\begin{equation*}
\K^{rc}\subset T_0+\mathrm{Span}\lt\{V_1,V_2,V_3\rt\}. 
\end{equation*}
This is straightforward because convex functions 
are Rank-$1$ convex. Thus $\K^{rc}\subset \mathrm{Conv}(K)\subset T_0+\mathrm{Span}\lt\{V_1,V_2,V_3\rt\}$. One general principle is, if $\VI:=\mathrm{Span}\lt\{V_1,V_2,V_3\rt\}$ does not contain enough Rank-$1$ directions, then $\K$ does not contain $\Tb_4$. This is the content of Lemma \ref{L13}. Because of Remark \ref{r:nodegenerateT4}, we only need to consider two cases: $\mathrm{dim}\lt(\VI\rt)=2$ and $\mathrm{dim}\lt(\VI\rt)=3$. The arguments to deal with the two cases are somewhat different and we will discuss each in turn.

\subsection{Case 1:  $\mathrm{dim}\lt(\VI\rt)=3$} 
\blue{An important observation is that if a linear isomorphism preserves Rank-$1$ matrices, then it preserves $\Tb_4$. This is the content of Lemma \ref{LL0}.} This fact allows us to transform the original set $\K$ into a simpler set $\Ub_{\K}^0$ given by
\begin{equation*}
	\Ub_{\KKI}^{0}:=\lt\{\lt(\begin{array}{cc} h_i & r_i \\ a(r_i) & h_i \\ h_i a(r_i) & \frac{h_i^2}{2}+F(r_i)\end{array}\rt): i=0,1,2,3\rt\},
\end{equation*}
where $h_i:=u_i-u_0$, $r_i:=v_i-v_0$ and \blue{the functions $a$ and $F$ are translations of the functions $\al$ and $\Fk$ satisfying the normalization $a(0)=F(0)=0$}.
By relatively straightforward arguments we can show that, denoting $\vec{h}=(h_1,h_2,h_3), \vec r=(r_1,r_2,r_3)$ and $\vec z=(a(r_1), a(r_2), a(r_3))$, if $\vec h\times\vec r=0$ or $\vec h\times\vec z=0$ then $\Ub_{\KKI}^{0}$ cannot contain a $\Tb_4$. So we can assume this is not the case. For $x,y\in \R^3$, let $\lt(x|y\rt)\in M^{3\times 2}$ denote the matrix whose columns are $x$ and $y$. A crucial observation is that if for some matrix $\AI\in M^{3\times 3}$ we can represent $\mathrm{Span}\{\Ub_{\KKI}^{0}\}$ in the form 
\begin{equation}
\label{intreq19}
\mathrm{Span}\{\Ub_{\KKI}^{0}\}=\lt\{\lt(z|\AI z\rt):z\in \R^3\rt\},
\end{equation}
then $M\in \mathrm{Span}\{\Ub_{\KKI}^{0}\}$ is Rank-$1$ if and only if $M=\lt(\zeta|\AI\zeta\rt)$ where $\zeta\in \R^3$ is an eigenvector of 
$\AI$. \blue{So if (\ref{intreq19}) holds, then the Rank-$1$ directions are contained in the eigenspaces of $\AI$, and thus, in the worst case, can form either a two-dimensional subspace and a line, or three distinct lines.  In either of these two cases, there are not enough Rank-$1$ directions to build three-dimensional $\Tb_4$ (see Lemma \ref{L13} (b); the above discussions are ideas of V. \v{S}ver\'{a}k communicated to the first author \cite{svper2}).} So the issue becomes to what extent we can write $\mathrm{Span}\{\Ub_{\KKI}^{0}\}$ in the form of (\ref{intreq19}). We can clearly find matrices $\AI_1, \AI_2\in M^{3\times 3}$ such that $\mathrm{Span}\{\Ub_{\KKI}^{0}\}=\lt\{\lt(\AI_1 z|\AI_2 z\rt):z\in \R^3\rt\}$. If either $\AI_1$ or $\AI_2$ is invertible then  $\mathrm{Span}\{\Ub_{\KKI}^{0}\}$ can be represented in the form of (\ref{intreq19}) and we are done (see Lemma \ref{L12}). Otherwise, letting $\lt(\AI_1|\AI_2\rt)\in M^{3\times 6}$ denote the matrix whose first three columns are the columns of $\AI_1$ and second three are the columns of $\AI_2$, we have two further cases to consider.

\subsubsection{The case $\mathrm{Rank}(\AI_1)=\mathrm{Rank}(\AI_2)=2$ and  $\mathrm{Rank}\lt(\lt(\AI_1|\AI_2\rt)\rt)=3$ (see Lemma \ref{L15}).}  In this case using the particular forms of $\AI_1$ and $\AI_2$ there exist
$\lm_1, \lm_2, \mu_1, \mu_2$ with $(\lm_1,\lm_2)\ne(\mu_1,\mu_2)$ such that 
\begin{equation*}
\mathrm{Span}\{\Ub_{\KKI}^{0}\}=\lt\{\lt(\begin{array}{cc} \vec h\cdot\vec\alpha & \vec r\cdot\vec\alpha \\ \vec z\cdot\vec\alpha & \vec h\cdot\vec\alpha \\ \lm_1(\vec h\cdot\vec\alpha) + \lm_2 (\vec z\cdot\vec\alpha) & \mu_1 (\vec r\cdot\vec\alpha) + \mu_2 (\vec h\cdot\vec\alpha)\end{array}\rt): \vec\alpha\in\R^3\rt\}.
\end{equation*}
\blue{Again the Rank-$1$ directions must satisfy $M_{ij}=0$ for all $i\ne j$. Similar to Case 1, a careful but straightforward analysis using the special structure of the three minors and the fact that $(\lm_1,\lm_2)\ne(\mu_1,\mu_2)$ shows that there are not enough Rank-$1$ directions in $\mathrm{Span}\{\Ub_{\KKI}^{0}\}$ to form three-dimensional $\Tb_4$.}

\subsubsection{The case  $\mathrm{Rank}\lt(\lt(\AI_1|\AI_2\rt)\rt)=2$}\label{subsec:ranktwo}
This turns out to be the hardest case. 
In this case using the particular forms of $\AI_1$ and $\AI_2$ there exist $\lm_1, \lm_2$ such that 
\begin{equation}
\label{eqint11}
h_i a(r_i)=\lm_1 h_i+\lm_2 a(r_i),\qd \frac{h_i^2}{2}+F(r_i)=\lm_1 r_i+\lm_2 h_i. 
\end{equation}
\blue{Since the third rows of the matrices in $\Ub_{\K}^0$ are linear combinations of the first two rows with the same multiplicity constants, it is not hard to show that it suffices to show the set
\begin{equation*}
\ti\Ub_{\KKI}^{0}:=\lt\{\lt(\begin{array}{cc}0 & 0\\ 0& 0\end{array}\rt), \lt(\begin{array}{cc}h_1& r_1\\ a(r_1)& h_1\end{array}\rt),\lt(\begin{array}{cc}h_2 & r_2\\ a(r_2)& h_2\end{array}\rt),\lt(\begin{array}{cc}h_3 & r_3\\ a(r_3)& h_3\end{array}\rt)\rt\}
\end{equation*}
does not contain a $\Tb_4$. The set $\ti\Ub_{\K}^0$ is a subset of $M^{2\times 2}$ and much more is known about $\Tb_4$ \rred{configurations} in $M^{2\times 2}$. In particular a result in \cite{sz} implies that, labeling the matrices in $\ti\Ub_{\KKI}^{0}$ by $\ti T_i$, if for some $i$, 
\begin{equation}\label{ep200}
\text{the set }\{\det(\ti T_i-\ti T_j)\}\text{ does not change sign for } j\ne i,
\end{equation}
then $\ti\Ub_{\KKI}^{0}$ does not contain a $\Tb_4$. So our goal is to establish (\ref{ep200}) for the set $\ti\Ub_{\KKI}^{0}$.}

\blue{Now \rred{comes} another \rred{important} idea. The set $\Ub_{\KKI}^0$ is defined with respect to the point $(u_0, v_0)$. However, a closer look at the whole process, one observes that there is no unique role played by $(u_0, v_0)$ and all previous arguments also apply to the set $\Ub_{\K}^{k}$ for $k=1,2,3$, where the set $\Ub_{\K}^k$ is the analog of $\Ub_{\K}^0$ but defined with respect to the point $(u_k,v_k)$, i.e.,
\begin{equation*}
\Ub_{\KKI}^{k}:=\lt\{\lt(\begin{array}{cc} h_i^k & r_i^k \\ a_k(r_i^k) & h_i^k \\ h_i^k a_k(r_i^k) & \frac{(h_i^k)^2}{2}+F_k(r_i^k)\end{array}\rt): i=0,1,2,3\rt\},
\end{equation*}
where $h_i^k:=u_i-u_k$, $r_i^k:=v_i-v_k$ and the functions $a_k$ and $F_k$ are translations of the functions $\al$ and $\Fk$ satisfying the normalization $a_k(0)=F_k(0)=0$. This observation allows us the extra power to assume all $(h_i^k,r_i^k)$ satisfies the system (\ref{eqint11}) with constants $\lm_1^k, \lm_2^k$ and this turns out to be crucial.}

To establish (\ref{ep200}) we assume without loss of generality $v_0< v_1< v_2< v_3$ (the case of equalities easily leads to a degenerate case). Let 
$D^k_i:=(h^k_i)^2- r^k_i a_{k}(r^k_i)$ and it is not hard to show $D^k_i=D^i_k$. Now we form the symmetric matrix 
\begin{equation*}
\mathcal{S}:=\lt(\begin{array}{cccc} 0 & D^0_1 & D^0_2 & D^0_3 \\
D^1_0 & 0 & D^1_2  & D^1_3 \\
D^2_0 & D^2_1 & 0  & D^2_3 \\
D^3_0 & D^3_1 & D^3_2 & 0 
 \end{array}\rt).
\end{equation*}
\rrred{Now (\ref{ep200}) reinterpreted for matrix $\mathcal{S}$ says that if $\Ub_{\KKI}^{0}$ contains a $\Tb_4$ then every row and column of $\mathcal{S}$ must change sign.} \blue{In Lemmas \ref{l20}-\ref{l22}, we establish some elementary properties about the structure of solutions to a system of the form (\ref{eqint11}). Using these properties and the fact
$0<r^0_1<r^0_2<r^0_3$, any attempt to fill out the entries of matrix $\mathcal{S}$ leads to a configuration \rrred{in which one row or column of $\mathcal{S}$ does not change sign and hence (\ref{ep200}) is satisfied} for some $i$ (see Lemma \ref{L23})}.

\subsection{Case 2:  $\mathrm{dim}\lt(\VI\rt)=2$} 

\rrred{By the assumption $\mathrm{dim}(\VI)=2$, we have} $\mathrm{dim}\lt(\mathrm{Span}\{\Ub_{\KKI}^{0} \}\rt)=2$. \rrred{Thus} there exist
$\gamma_1,\gamma_2$, $\lm_1$, $\lm_2$ and $\mu_1$, $\mu_2$ such that $r_i=\gamma_1 h_i + \gamma_2 a(r_i)$, 
$h_i a(r_i)=\lm_1 h_i + \lm_2 a(r_i)$ and $\frac{h_i^2}{2}+F(r_i)=\mu_1 h_i + \mu_2 a(r_i)$. Therefore
\begin{equation*}
	\mathrm{Span}\{\Ub_{\KKI}^{0}\}=\lt\{\OI(s,t):=\lt(\begin{array}{cc} s & \gamma_1 s+\gamma_2 t \\ t & s \\ \lm_1 s+\lm_2 t & \mu_1 s + \mu_2 t \end{array}\rt): s,t\in\R\rt\}.
\end{equation*}
The Rank-$1$ directions required 
to build the $\Tb_4$ are contained in this subspace and must satisfy $M_{12}=M_{13}=M_{23}=0$, where $M_{ij}(P)$ denotes the $2\times 2$ minor of matrix $P\in M^{3\times 2}$ which is 
comprised of the $i$-th and $j$-th rows. So
\begin{equation*}
	M_{12}\lt(\OI(s,t)\rt)=s^2-\gamma_1 st-\gamma_2 t^2. 
\end{equation*}
If the discriminant $\gamma_1^2+4\gamma_2\leq 0$ then clearly there are not enough Rank-$1$ directions in $\mathrm{Span}\{\Ub_{\KKI}^{0}\}$ to build $\Tb_4$. So we must have $\gamma_1^2+4\gamma_2>0$ and hence $s^2-\gamma_1 st+\gamma_2 t^2=\lt(s-k t\rt)\lt(s-l t\rt)$ for some $k\ne l$. Thus the two possible Rank-$1$ directions are $\OI(kt, t)$ and $\OI(lt, t)$. In order for these two candidates to be Rank-$1$ directions, they must further satisfy $M_{13}=M_{23}=0$. Using the special structures of the three minors, one can show that if $\OI(kt, t)$ and $\OI(lt, t)$ are Rank-$1$ directions, then the coefficients $\gamma_1,\gamma_2$, $\lm_1$, $\lm_2$ and $\mu_1$, $\mu_2$ must satisfy a system of equations. As a consequence of simple algebraic manipulations, it turns out that $(h_i, r_i)$ must satisfy the system \eqref{eqint11}, and thus the arguments in Subsection \ref{subsec:ranktwo} allow us to conclude the proof.

\section{Preliminaries}


%
%

%
%

%
%

%
%

\blue{In what follows, we make the following convention. Given a set $\K:=\{T_i\}_{i=1}^N\subset M^{m\times n}$, we say that $\K$ does not \emph{contain} a $\Tb_N$ configuration if \emph{any} ordering of the elements in $\K$ cannot form a $\Tb_N$ configuration.} We first recall the following convenient result which is \blue{an immediate consequence of} Proposition 1 in \cite{sz} and characterizes $\Tb_N$ configurations in $M^{2\times 2}$. 
\begin{a5}[\cite{sz}]
\label{rmk2}
	\blue{Given a set $\{T_i\}_{i=1}^N\subset M^{2\times 2}$, a necessary condition for the set to \blue{contain} a $\Tb_N$ configuration is that, for every $i$, the set $\{\det(T_i-T_j):j\ne i\}$ changes sign.}
\end{a5}

\begin{a1}
	\label{L1}
	Given $\KKI:=\{T_1,\dots, T_N\}\subset M^{m\times n}$, let 
		\begin{equation}
		\label{ggeqa1}
		\green{V_{k}:=T_{k}-T_1,\qd k=2,3,\dots N,}
		\end{equation}
	and denote $\bblue{\VI}:=\mathrm{Span}\lt\{V_2,V_3,\dots, V_{N}\rt\}$. Then 
	\begin{equation*}
	\KKI^{rc}\subset \green{T_1}+\bblue{\VI}.
	\end{equation*}
\end{a1}
\begin{proof}
	\blue{Since convex functions 
	are Rank-$1$ convex, it follows that $\K^{rc}\subset \mathrm{Conv}(\K)\subset T_1+ \VI$.}
\end{proof}

\begin{a1}
\label{LL0}
Let $\VI\subset M^{m\times n}$ be a subspace and
$L:\VI\rightarrow \WI\subset M^{\blue{p\times q}}$ be a linear isomorphism with the property that 
\begin{equation}
\label{eqa30}
\mathrm{Rank}(A)=1\rblue{\Longleftrightarrow} \mathrm{Rank}\lt(L(A)\rt)=1.
\end{equation}
Then
\begin{equation*}
\lt\{T_1,\dots,T_N\rt\}\subset \VI\text{ forms a }\mathbb{T}_N\;
\rblue{\Longleftrightarrow} \lt\{L(T_1), \dots, L(T_N)\rt\}\subset \WI\text{ forms a }\mathbb{T}_N. 
\end{equation*}
\end{a1}

\begin{proof}
\rblue{We only need to establish the forward implication, as the reverse one follows the same lines by noting that $L^{-1}$ is a linear isomorphism satisfying (\ref{eqa30}) provided that $L$ satisfies (\ref{eqa30}).} 

Assume $\KKI:=\lt\{T_1,\dots,T_N\rt\}\subset \VI$ forms a $\mathbb{T}_N$, then there 
exist $P\in M^{m\times n}$, Rank-$1$ matrices $C_i\in M^{m\times n}$ and scalars $\kappa_i>1$ such that 
(\ref{eq1}) and \blue{(\ref{eq3.5})} hold true. \blue{Defining $V_k$'s as in (\ref{ggeqa1}), it is clear that $V_k\in \VI$ and thus it follows from Lemma \ref{L1} that} 
\begin{equation}
\label{eq529aabb1}
\blue{\KKI^{rc}\subset T_1+\mathrm{Span}\lt\{V_2,V_3,\dots, V_{N}\rt\}\subset \VI}.
\end{equation}
\rred{Let \blue{the} matrices $\lt\{P_i\rt\}$ be defined by 
\begin{equation*}
\blue{P_i=P+C_1+\dots+C_{i-1},}
\end{equation*}
where \blue{$P$ and $C_i$ are as in Definition \ref{def1} and} the index $i$ is counted modulo $N$. Then \blue{as} shown in the paragraph after Definition 7 of \cite{kms}, we have that each $P_i\in \KKI^{rc}$. In particular, as
\begin{equation}\label{eq529aabb2}
\blue{C_i=P_{i+1}-P_{i}},
\end{equation}
we have}
\begin{equation}
\label{eq529aabb3}
\rred{C_i \overset{(\ref{eq529aabb1}), (\ref{eq529aabb2})}{\in}  \blue{\VI}.}
\end{equation}
Now by (\ref{eqa30}) we have that $L(C_i)$ is Rank-$1$ and by linearity of $L$ we have that 
$\lt\{L(T_1), \dots, L(T_N)\rt\}$ satisfies (\ref{eq1}) for $L(P)$, $L(C_i)$, $\kappa_i$ for $i=1,\dots,N$. \rblue{Further, $\lt\{L(T_1), \dots, L(T_N)\rt\}$ has no Rank-$1$ connections as a result of (\ref{eqa30}) and the fact that $\KKI$ contains no Rank-$1$ connections. Thus $\lt\{L(T_1), \dots, L(T_N)\rt\}$ forms a $\mathbb{T}_N$.}
\end{proof}

%
%

%
%

\blue{For the rest of this paper, we will focus on $\Tb_4$ configurations in the set $\KKI_1$ defined in (\ref{ep100}) under the assumption that the function $\al$ is monotonic increasing and strictly convex, i.e., $\al'>0$ and $\al''>0$, unless otherwise specified. Given a set $\KKI$ of four points in $\KKI_1$, for technical reasons, it is more convenient for most of the time to label the four points as $T_i=P(u_i,v_i)$ for $i=0,1,2,3$, where recall that the mapping $P:\R^2\rightarrow \KKI_1$ is defined in (\ref{eq4}), and thus
\begin{equation}\label{ep7}
\KKI=\{P(u_0,v_0),P(u_1,v_1),P(u_2,v_2), P(u_3,v_3)\}.
\end{equation}	
We denote by 
\begin{equation}\label{eqp1}
h_i=u_i-u_0,\qd r_i=v_i-v_0, 
\end{equation}	
and $\vec h=(h_1,h_2,h_3), \vec r=(r_1,r_2,r_3)$. It should be pointed out that all the results in the remaining of this paper do not rely on any particular ordering of the four points. We first make some simplifications.}

\begin{a1}
\label{L5}
Given $\KKI$ as in (\ref{ep7}), define 
$V_i:=P(u_i,v_i)-P(u_0,v_0)$. There exists an invertible matrix $B\in M^{3\times 3}$ such that 
\begin{equation}
\label{eq31}
B V_i=\lt(\begin{array}{cc} h_i & r_i \\
\al(v_0+r_i)-\al(v_0) & h_i \\
h_i(\al(v_0+r_i)-\al(v_0)) & \frac{h_i^2}{2}+\Fk(v_0+r_i)-\Fk(v_0)-\al(v_0)r_i \end{array}\rt).
\end{equation}
\end{a1}
\begin{proof} 
Using (\ref{eq4}) we write
\begin{eqnarray*}
V_i&=&\lt(\begin{matrix}
h_i & r_i\\
\al(v_0+r_i)-\al(v_0)  & h_i\\
(u_0+h_i)\al(v_0+r_i)-u_0\al(v_0)  & u_0 h_i+\frac{h_i^2}{2}+\Fk(v_0+r_i)-\Fk(v_0)
\end{matrix}\rt).
\end{eqnarray*}
Multiplying the second row by $u_0$ and subtracting it from the third row we obtain 
\begin{equation*}
\hat{V}_i=\lt(\begin{matrix}
h_i & r_i\\
\al(v_0+r_i)-\al(v_0)  & h_i\\
 h_i\al(v_0+r_i)&  \frac{h_i^2}{2}+\Fk(v_0+r_i)-\Fk(v_0)
\end{matrix}\rt).
\end{equation*}
Multiplying the first row by $a(v_0)$ and subtracting it from the third row in $\hat V_i$ we obtain 
\begin{equation}
\label{eq54}
\hat{\hat{V}}_i=\lt(\begin{array}{cc} h_i & r_i \\
\al(v_0+r_i)-\al(v_0) & h_i \\
h_i(\al(v_0+r_i)-\al(v_0)) & \frac{h_i^2}{2}+\Fk(v_0+r_i)-\Fk(v_0)-\al(v_0)r_i \end{array}\rt).\nn
\end{equation}
This establishes (\ref{eq31}). 
\end{proof}

%
%

To simplify notation, for a fixed $v\in\R$, define
\begin{equation}
\label{eq6}
a_{\blue{v}}(t) := \al(v+t)-\al(v), \qd F_{\blue{v}}(t):=\Fk(v+t)-\Fk(v)-\al(v)t.
\end{equation}
Since $\al'>0, \al''>0$ and $\Fk'=\al$, it is clear that
\begin{equation}
\label{eq7}
a_{v}(0)=0, \qd a_{v}'(t)>0, \qd a_{v}''(t)>0
\end{equation}
and
\begin{equation}
\label{eq8}
F_{v}'(t)=a_{v}(t),\qd F_{v}''(t)=a_{v}'(t)>0,\qd F_{v}(0)=F_{v}'(0)=0.
\end{equation}
Further, given $h,r\in \R$, define 
\begin{equation}
\label{eq400}
\QI_{\blue{v}}(h,r):=\lt(\begin{array}{cc} h & r \\
a_{\blue{v}}(r)  & h \\
ha_{\blue{v}}(r) & \frac{h^2}{2}+F_{\blue{v}}(r) \end{array}\rt).
\end{equation}
For $\KKI$ given in (\ref{ep7}), \blue{we define the associated set $\Ub_{\KKI}^{0}$ with respect to the point $P(u_0,v_0)$ by}
\begin{equation}
\label{eqe25}
\blue{\Ub_{\KKI}^{\blue{0}}}:=\lt\{ \QI_{\blue{v_0}}(0,0) , \QI_{\blue{v_0}}(h_1,r_1), \QI_{\blue{v_0}}(h_2,r_2),\QI_{\blue{v_0}}(h_3,r_3)\rt\},
\end{equation}
where $h_i, r_i$ are defined in (\ref{eqp1}). We will need the following fundamental result.
\begin{a1}
\label{L1.2}
If $\KKI$ (given in (\ref{ep7})) \blue{contains} a $\Tb_4$, then $\Ub_{\KKI}^{\blue{0}}$ also
\blue{contains} a $\Tb_4$. 
\end{a1}

\begin{proof}
	\blue{Without loss of generality, we may assume that the ordering $\{T_0, T_1, T_2, T_3\}$ forms a $\Tb_4$. Denoting $T_i:=P(u_i,v_i)$ and $V_i=T_i-T_0$, it is clear that $\{0, V_1, V_2, V_3\}\subset M^{3\times 2}$ forms a $\Tb_4$. Now we define $\VI:=\mathrm{Span}\{V_1, V_2, V_3\}$ and the linear mapping $L: \VI\rightarrow M^{3\times 2}$ by $L(X)=BX$, where $B\in M^{3\times 3}$ is the invertible matrix found in Lemma \ref{L5}. Since the mapping $L$ corresponds to row operations, it is clearly a linear isomorphism satisfying (\ref{eqa30}). The lemma follows from Lemmas \ref{LL0} and \ref{L5}.}
\end{proof}

\section{Non-existence of $\Tb_4$ in some special cases}

In this section, given $\KKI$ as in (\ref{ep7}), we show that the four points cannot \blue{contain} a $\Tb_4$ if the vectors $\vec h$ and $\vec r$ defined in (\ref{eqp1}) satisfy certain special relations. By Lemma \ref{L1.2}, it is sufficient to show that the set $\Ub_{\KKI}^{0}$ defined in (\ref{eqe25}) cannot \blue{contain} a $\Tb_4$. To simplify notation, when there is no risk of confusion, we omit the dependence of the mapping $\QI_{v}$ and the functions $a_{v}, F_{v}$ on $v$. Let
\begin{equation*}
\Lambda_R:=\lt\{A\in M^{3\times 2}:\mathrm{Rank}(A)=1\rt\},
\end{equation*}
i.e., the cone of all Rank-$1$ matrices in $M^{3\times 2}$. 

\begin{a1}\label{L13}
	Let $\Ub_{\KKI}^{0}$ be defined by (\ref{eqe25}).
	\begin{enumerate}[label=(\alph*)]
		\item \red{If $\mathrm{dim}\lt(\mathrm{Span}\{\Ub_{\KKI}^{0}\}\rt)=2$ and $\Lambda_R\cap \mathrm{Span}\{\Ub_{\KKI}^{0}\}$ consists of a single line then $\Ub_{\KKI}^{0}$ cannot \blue{contain} a $\Tb_4$.}
	 	\item If $\mathrm{dim}\lt(\mathrm{Span}\{\Ub_{\KKI}^{0}\}\rt)=3$ and $\Lambda_R\cap \mathrm{Span}\{\Ub_{\KKI}^{0}\}$
		either consists of at most three distinct lines or a two-dimensional plane and a line, then $\Ub_{\KKI}^{0}$ cannot \blue{contain} a $\Tb_4$. 
	\end{enumerate}
\end{a1}

\begin{proof} 
	\blue{The proof of (a) is trivial. We focus on (b) and assume $\mathrm{dim}\lt(\mathrm{Span}\{\Ub_{\KKI}^{0}\}\rt)=3$.} Suppose  $\Lambda_R\cap \mathrm{Span}\{\Ub_{\KKI}^{0}\}$ consists of three 
	distinct lines and \blue{without loss of generality assume that $\Ub_{\KKI}^{0}$ with the given ordering forms a $\Tb_4$, then there exist $C_i\in\Lambda_R, i=0,1,2,3, P\in M^{3\times 2}, \kappa_i>1$ such that (\ref{eq1}) and (\ref{eq3.5}) hold true. By Lemma \ref{L1} \rred{and (\ref{eq529aabb3}) we have}  $C_i\in \Lambda_R\cap \mathrm{Span}\{\Ub_{\KKI}^{0}\}$}. Thus, for some $i_0\ne i_1\in \lt\{0, 1,2,3\rt\}$, there exists $\lm\blue{\ne 0}$ such that $C_{i_1}=\lm C_{i_0}$. Let $i_2,i_3$ be such that $\lt\{i_2,i_3\rt\}=\lt\{0, 1,2,3\rt\}\backslash \lt\{i_0,i_1\rt\}$. Equation (\ref{eq3.5}) then becomes 
	\begin{equation}
	\label{eq172}
	(1\blue{+}\lm)C_{i_0}+C_{i_2}+C_{i_3}=0.  
	\end{equation}
	So the matrices $C_{i_0}, C_{i_2}, C_{i_3}$ are linearly dependent and their span forms a \blue{subspace $\VI$ of dimension at most two.
	It follows from (\ref{eq1}) that 
	\begin{equation}
	\label{eq174}
	\Ub_{\KKI}^{0}\subset P+\VI.
	\end{equation}
	Now since $\QI(0,0)=0\in P+\VI$, it is clear that $P+\VI$ is a subspace of dimension at most two, and this contradicts our assumption that $\mathrm{dim}\lt(\mathrm{Span}\{\Ub_{\KKI}^{0}\}\rt)=3$.} 
	
	Next suppose   $\Lambda_R\cap \mathrm{Span}\{\Ub_{\KKI}^{0}\}$  consists of a two-dimensional plane $\WI$ and a single line $\LI\nsubseteq\WI$ and \blue{again assume $\Ub_{\KKI}^{0}$ with the given ordering forms a $\Tb_4$. Let $C_i, P, \kappa_i$ be as above. If $C_i\in \WI$ for all $i$, then similar to (\ref{eq174}) we have $\Ub_{\KKI}^{0}\subset P+\WI$ and thus $\mathrm{dim}\lt(\mathrm{Span}\{\Ub_{\KKI}^{0}\}\rt)\leq 2$, which is a contradiction.}
	Let $i_0\in \lt\{0,1,2,3\rt\}$ be such that $C_{i_0}\in\LI$. If $C_i\in\WI$ for all $i\ne i_0$, then (\ref{eq3.5}) implies $C_{i_0}=-\sum_{i\ne i_0} C_{i}\in \WI$, which is a contradiction. \blue{So there exists $i_1\in\{0,1,2,3\}\setminus\{i_0\}$ such that $C_{i_1}\in\LI$ and thus $C_{i_1}=\lm C_{i_0}$ for some $\lm\ne 0$}. \rred{Thus equation (\ref{eq172}) must be satisfied and arguing exactly as in the last paragraph this contradicts the assumption that $\mathrm{dim}\lt(\mathrm{Span}\{\Ub_{\KKI}^{0}\}\rt)=3$. This completes the proof.}
\end{proof}

%
%

\blue{For the rest of this paper, besides the notations $\vec h=(h_1,h_2,h_3), \vec r=(r_1,r_2,r_3)$, we will further use
\begin{equation}\label{eqp2}
\vec z:=\lt(a(r_1), a(r_2), a(r_3)\rt), \qd\vec y:=\lt(h_1 a(r_1), h_2a(r_2), h_3a(r_3)\rt),
\end{equation}
and
\begin{equation}\label{eqp3}
\vec w:=\lt(\frac{h_1^2}{2}+F(r_1),\frac{h_2^2}{2}+F(r_2),\frac{h_3^2}{2}+F(r_3)\rt).
\end{equation}
And we will use $(\hat\cdot)$ to denote two-dimensional vectors.}

\begin{a1}
\label{L1.3}
Let $\Ub_{\KKI}^{0}$ be defined by (\ref{eqe25}). If $\vec h=0$ or $\vec r=0$, then $\Ub_{\KKI}^{0}$ cannot \blue{contain} a $\Tb_4$. 
\end{a1}
\begin{proof}
\green{\em Case 1. \rm We start by considering the case $\vec r=0$.}

\green{\em Proof of Case 1.\rm} First note that for $i_1\ne i_2\in \lt\{1,2,3\rt\}$ we have $h_{i_1}\not=h_{i_2}$ \blue{since 
otherwise $\ca{\Ub_{\KKI}^{0}}\leq 3$}. For the same reason we have $h_i\not=0$ for any $i=1,2,3$. Now 
$\det\lt(\begin{array}{cc} h_1 & h_2 \\ h_1^2 & h_2^2 \end{array}\rt)=h_1 h_2(h_2-h_1)\not=0$ and thus 
$\lt(\begin{array}{c} h_1 \\ h_1^2/2 \end{array}\rt)$ and $\lt(\begin{array}{c} h_2 \\ h_2^2/2 \end{array}\rt)$ are linearly independent. Let 
$\hat h:=\lt(h_1,h_2\rt)$ and $\hat w=\lt(\frac{h_1^2}{2},\frac{h_2^2}{2}\rt)$. \green{Since \blue{$\lt(\begin{array}{c} h_3 \\ h_3^2/2 \end{array}\rt)\in \mathrm{Span}\lt\{\lt(\begin{array}{c} h_1 \\ h_1^2/2 \end{array}\rt),\lt(\begin{array}{c} h_2 \\ h_2^2/2 \end{array}\rt)\rt\}$},} we have
\begin{eqnarray*}
\mathrm{Span}\{\Ub_{\KKI}^{0}\}&\overset{(\ref{eqe25}), (\ref{eq400})}{=}&\mathrm{Span}\lt\{\lt(\begin{array}{cc} h_1 & 0 \\ 0 & h_1 \\ 0 & \blue{\frac{h_1^2}{2}}  \end{array}\rt),  \lt(\begin{array}{cc} h_2 & 0 \\ 0 & h_2 \\ 0 & \blue{\frac{h_2^2}{2}}  \end{array}\rt) \rt\}\nn\\
&=&\lt\{\lt(\begin{array}{cc} \hat h\cdot \hat\alpha & 0 \\ 0 & \hat h\cdot \hat\alpha \\ 0 & \hat w\cdot \hat\alpha  \end{array}\rt):\hat\alpha\in \R^2\rt\}.
\end{eqnarray*}
Note that $\mathrm{Rank}\lt(\begin{array}{cc} \hat h\cdot \hat\alpha & 0 \\ 0 & \hat h\cdot \hat\alpha \\ 0 & \hat w\cdot \hat\alpha  \end{array}\rt)=1$ if and only if $\hat h\cdot \hat\alpha=0$. So there is only one Rank-$1$ line inside $\mathrm{Span}\{\Ub_{\KKI}^{0}\}$ and thus \blue{Lemma \ref{L13} (a) completes the proof in Case 1.}\nl

\green{\em Case 2. \rm We consider the case where $\vec h=0$ and $\mathrm{dim}\lt(\mathrm{Span}\{\Ub_{\KKI}^{0}\}\rt)=2$}. 

\green{\em Proof of Case 2. \rm } Now we have
\begin{equation*}
\QI(0,r_i)\overset{(\ref{eq400})}{=}\lt(\begin{matrix}
0 & r_i\\
a(r_i)  & 0\\
0 &  F(r_i)
\end{matrix}\rt).
\end{equation*}
Without loss of generality, assume that
\begin{eqnarray}\label{eq12}
\mathrm{Span}\{\Ub_{\KKI}^{0}\}&=&\mathrm{Span}\lt\{\lt(\begin{array}{cc} 0 & r_1 \\ a(r_1) & 0 \\ 0 & F(r_1)  \end{array}\rt),  \lt(\begin{array}{cc} 0 & r_2 \\ a(r_2) & 0 \\ 0 & F(r_2)  \end{array}\rt) \rt\}\nn\\
&=&\lt\{\lt(\begin{array}{cc} 0 & \hat r\cdot \hat\alpha \\ \hat z\cdot\hat\alpha & 0 \\ 0 & \hat w\cdot \hat\alpha  \end{array}\rt):\hat\alpha\in \R^2\rt\}.
\end{eqnarray}
We claim that 
\begin{equation}
\label{aaeqbb3}
\green{\hat r\text{ and }\hat w\text{ are linearly independent}}. 
\end{equation}
Suppose not, then there exists \rblue{$\lm\in\R$} such that 
\begin{equation*}
F(r_i)=\lm r_i \qd\text{ for } i=1,2,
\end{equation*} 
and therefore $r_i$ is a root of $g(t):=F(t)-\lm t$. Note that $g'(t)=a(t)-\lm$ \blue{and $g''(t)=a'(t)>0$ by (\ref{eq7}), and thus the function $g$ is strictly convex and has at most two roots.} It is clear that $g(0)=0$ using (\ref{eq8}), and thus \blue{$r_1=0$ or $r_2=0$} which \rrred{ as in Case 1 implies $\ca{\Ub_{\KKI}^{0}}\leq 3$} \bblue{and}  is a contradiction. So (\ref{aaeqbb3}) is established. Note that there are only two non-trivial minors in $\mathrm{Span}\{\Ub_{\KKI}^{0}\}$, namely,
\begin{equation*}
M_1 = (\hat r\cdot\hat\alpha)(\hat z\cdot\hat\alpha) \qd\text{and}\qd M_2=(\hat z\cdot\hat\alpha)(\hat w\cdot\hat\alpha).
\end{equation*}
So the Rank-$1$ directions must satisfy $M_1=M_2=0$. This requires either 
\begin{equation}\label{eqp17}
\hat z\cdot\hat\alpha =0
\end{equation}
or 
\begin{equation}\label{eqp18}
\hat r\cdot \hat\alpha=0 \qd\text{and}\qd \hat w\cdot\hat\alpha=0.
\end{equation} 
In the latter case, because of (\ref{aaeqbb3}), there is no Rank-$1$ direction. \blue{Clearly (recalling \rred{(\ref{eqp2})}) $\hat z\ne0\in\R^2$}, hence there is only one Rank-$1$ direction in $\mathrm{Span}\{\Ub_{\KKI}^{0}\}$ from the equation (\ref{eqp17}). \blue{We appeal to Lemma \ref{L13} (a) again to complete Case 2.}\nl

\green{\em Case 3. \rm We consider the case where $\vec h=0$ and $\mathrm{dim}\lt(\mathrm{Span}\{ \Ub_{\KKI}^{0}\}\rt)=3$}.

\green{\em Proof of Case 3. \rm \blue{Following exactly the same lines as in Case 2, we have an analogous expression for $\mathrm{Span}\{\Ub_{\KKI}^{0}\}$ as in (\ref{eq12}) with two-dimensional vectors replaced by three-dimensional vectors}, and $\vec r$ and $\vec w$ are linearly independent. As in (\ref{eqp17}) and (\ref{eqp18}), the Rank-$1$ directions in $\mathrm{Span}\{\Ub_{\KKI}^{0}\}$ must satisfy}
\begin{equation*}
\vec z\cdot\vec\alpha =0
\end{equation*}
or
\begin{equation*}
\vec r\cdot \vec\alpha=0 \qd\text{and}\qd \vec w\cdot\vec\alpha=0.
\end{equation*} 
In the first case, the Rank-$1$ directions form a two-dimensional plane. In the second case, as $\vec r$ and $\vec w$ are linearly independent, there is only one Rank-$1$ line. \rred{So the entire set of Rank-$1$ directions in $\mathrm{Span}\{\Ub_{\KKI}^{0}\}$ is the union of a two-dimensional plane and a line, and thus} \blue{we apply Lemma \ref{L13} (b) to finish the proof.}
\end{proof}


\begin{a1}
	\label{L1.5}
	Let $\Ub_{\KKI}^{0}$ be defined by (\ref{eqe25}). Recalling (\ref{eqp2}), if $\vec h\times \vec r=0$ or $\vec h\times\vec z=0$, \rred{then}
\begin{equation}
\label{eq529aabb4}
\rred{\Ub_{\KKI}^{0}\text{ cannot \blue{contain} a }\Tb_4. }
\end{equation}
\end{a1}

\begin{proof}
\rred{\em Step 1. \rm We will show (\ref{eq529aabb4}) under the assumption $\vec h\times\vec z=0$.}
	
\rred{\em Proof of Step 1. \rm We} may assume that $\vec h\ne 0$ and $\vec r\ne 0$ by Lemma \ref{L1.3} and hence $\vec z\overset{\rrred{(\ref{eqp2})}}{\ne} 0$. So there exists some $\lm\ne 0$ such that 
\begin{equation}
\label{eq529aabb5}
\rred{\vec z=\lm \vec h. }
\end{equation}
Thus 
	\begin{equation}
	\label{eqe71}
	\QI(h_i,r_i)\overset{\rred{(\ref{eq400})}}{=}\lt(\begin{array}{cc} h_i & r_i \\ \lm h_i & h_i \\ \lm h_i^2 & \frac{h_i^2}{2}+F(r_i)  \end{array}\rt)\qd\text{ for }i=1,2,3.
	\end{equation}
	
	First assume that $\dim\lt(\mathrm{Span}\{\Ub_{\KKI}^{0}\}\rt)=2$. Without loss of generality assume that $\QI(h_1,r_1)$ and $\QI(h_2,r_2)$ are linearly independent and thus (recalling (\ref{eqp2}) and (\ref{eqp3}))
	\begin{equation}\label{eq17}
	\mathrm{Span}\{\Ub_{\KKI}^{0}\}=\lt\{\lt(\begin{array}{cc} \hat h\cdot\hat \alpha & \hat r\cdot \hat\alpha \\ \lm \hat h\cdot\hat\alpha & \hat h\cdot \hat\alpha \\ \lm \hat p\cdot \hat\alpha & \hat w\cdot \hat\alpha  \end{array}\rt):\hat\alpha\in \R^2\rt\},
	\end{equation}
	where $\hat p=(h_1^2,h_2^2)$. If $\hat r\times\hat h=0$, then $\hat h=\mu\hat r$ for some $\mu\ne 0$ and $r_1, r_2$ are solutions of $a(t)\overset{\rred{(\ref{eq529aabb5}), (\ref{eqp2})}}{=}\lm\mu t$. However, \rred{as we have seen before }since $a$ is strictly convex, the equation has at most one non-trivial solution. If $r_i=0$ for some $i$, then $h_i=\mu r_i=0$; or if $r_1=r_2$, we have $h_1=h_2$. In both cases \rred{from (\ref{eqe71})} we have $\ca{\Ub_{\KKI}^0}\leq 3$. Similar arguments using the convexity of the square function show that $\ca{\Ub_{\KKI}^0}\leq 3$ if $\hat p\times \hat h=0$. So we can assume that
	\begin{equation}\label{ep1}
	\hat r\times\hat h\ne 0, \qd \hat p\times \hat h\ne 0.
	\end{equation}
	Note that the three minors in $\mathrm{Span}\{\Ub_{\KKI}^{0}\}$ are
	\begin{equation}\label{ep2}
	M_1 = (\hat h\cdot \hat\alpha)^2-\lm(\hat h\cdot\hat\alpha)(\hat r\cdot\hat\alpha),
	\end{equation}
	\begin{equation}
	M_2 = (\hat h\cdot\hat\alpha)(\hat w\cdot\hat\alpha)-\lm(\hat r\cdot\hat\alpha)(\hat p\cdot\hat\alpha),
	\end{equation}
	and
	\begin{equation}\label{ep3}
	M_3 = \lm(\hat h\cdot\hat\alpha)(\hat w\cdot\hat\alpha)-\lm(\hat h\cdot\hat\alpha)(\hat p\cdot\hat\alpha). 
	\end{equation}
	The Rank-$1$ directions in $\mathrm{Span}\{\Ub_{\KKI}^{0}\}$ must satisfy $M_1=M_2=M_3=0$.  From $M_1=0$, we need $\hat h\cdot\hat\alpha=0$ or $\hat h\cdot\hat\alpha=\lm \hat r\cdot\hat\alpha$. When $\hat h\cdot \hat\alpha=0$, it follows from $M_2=0$ that $\hat r\cdot\hat\alpha=0$ or $\hat p\cdot\hat\alpha=0$. Recall that we have (\ref{ep1}). Hence in this case we always have $\hat\alpha=0$ and thus there is no Rank-$1$ direction. When $\hat h\cdot\hat\alpha=\lm \hat r\cdot\hat\alpha$, we have $(\hat h-\lm \hat r)\cdot\hat\alpha=0$. By (\ref{ep1}) we know $\hat h-\lm \hat r\ne 0$, and hence there is at most one Rank-$1$ direction. Putting the above together, when $\dim\lt(\mathrm{Span}\{\Ub_{\KKI}^{0}\}\rt)=2$, there is at most one Rank-$1$ direction in $\mathrm{Span}\{\Ub_{\KKI}^{0}\}$ and thus Lemma \ref{L13} (a) applies.
	
	Now we assume that $\dim\lt(\mathrm{Span}\{\Ub_{\KKI}^{0}\}\rt)=3$. Then the expressions (\ref{eq17}) and (\ref{ep2})-(\ref{ep3}) still hold with two-dimensional vectors replaced by three-dimensional vectors. Following exactly the same lines \rred{of argument} as above, we may assume 
	\begin{equation}\label{ep4}
	\vec r\times \vec h\ne 0, \qd\vec p\times\vec h\ne 0.
	\end{equation}
	The Rank-$1$ directions still satisfy $M_1=M_2=M_3=0$. From $M_1=0$, we need $\vec h\cdot\vec\alpha=0$ or $\vec h\cdot\vec\alpha=\lm \vec r\cdot\vec\alpha$. When $\vec h\cdot \vec\alpha=0$, it follows from $M_2=0$ that $\vec r\cdot\vec\alpha=0$ or $\vec p\cdot\vec\alpha=0$. Because of (\ref{ep4}), there are at most two Rank-$1$ directions in this case. When $\vec h\cdot\vec\alpha=\lm \vec r\cdot\vec\alpha$, the set of Rank-$1$ directions satisfies $(\vec h-\lm \vec r)\cdot\vec \alpha=0$, and forms at most a two-dimensional plane thanks to (\ref{ep4}). Note that the Rank-$1$ direction determined by $\vec h\cdot\vec \alpha=0$ and $\vec r\cdot\vec\alpha=0$ is contained in this plane. Thus when $\dim\lt(\mathrm{Span}\{\Ub_{\KKI}^{0}\}\rt)=3$, the Rank-$1$ directions in $\mathrm{Span}\{\Ub_{\KKI}^{0}\}$ \rred{are contained in the union of }a line and at most a two-dimensional plane. This allows us to use Lemma \ref{L13} (b) to conclude the proof of Step 1. \nl

\em Step 2. \rm \rred{We will show (\ref{eq529aabb4}) under the assumption $\vec h\times\vec r=0$. }

\em Proof of Step 2. \rm  There exists some $\lm\ne 0$ such that 
\begin{equation}
\label{eq529aabb6}
\rred{\vec r=\lm \vec h. }
\end{equation}
Thus 
\begin{equation}
\label{eq529aabb8}
\QI(h_i,r_i)\overset{\rred{(\ref{eq400})}}{=}\lt(\begin{array}{cc} h_i & \lm h_i \\ a(r_i) & h_i \\ h_i a(r_i) & \frac{h_i^2}{2}+F(r_i)  \end{array}\rt)\qd\text{ for }i=1,2,3.
\end{equation}
	
	First assume that $\dim\lt(\mathrm{Span}\{\Ub_{\KKI}^{0}\}\rt)=2$. Again assume \rred{without loss of generality} that $\QI(h_1,r_1)$ and $\QI(h_2,r_2)$ are linearly independent and we obtain (recalling (\ref{eqp2}) and (\ref{eqp3}))
\begin{equation}
\label{ss43}
\mathrm{Span}\{\Ub_{\KKI}^{0}\}=\lt\{\lt(\begin{array}{cc} \hat h\cdot \hat\alpha & \lm \hat h\cdot \hat\alpha \\ \hat z\cdot\hat\alpha & \hat h\cdot \hat\alpha \\ \hat y\cdot \hat\alpha & \hat w\cdot \hat\alpha  \end{array}\rt):\hat\alpha\in \R^2\rt\}.
\end{equation}
Similar to the arguments in Step 1, \rred{we claim that 
\begin{equation}
\label{eq529aabb7}
\rrred{\hat z\times \hat h=0 \Longrightarrow \ca{\Ub_{\KKI}^{0}}\leq 3.} 
\end{equation}
\blue{Indeed} if \bblue{$\hat z\times \hat h=0$} we have $\hat{z}=\mu \hat{h}$ for some $\mu\not=0$, so $\hat{z}\overset{\rrred{(\ref{eq529aabb6})}}{=}\frac{\mu}{\lm}\hat{r}$ and by convexity of $a$ either this implies $r_i=0$ for some $i$ or 
$r_{i_0}=r_{i_1}$ for some $i_0\ne i_1$. In either case by (\ref{eq529aabb6}) and (\ref{eq529aabb8}), \rrred{ we have that (\ref{eq529aabb7}) follows}. }

\rrred{In a very similar way, we claim that} 
\begin{equation}
\label{ss44}
\rrred{\hat w\times \hat h=0 \Longrightarrow \ca{\Ub_{\KKI}^{0}}\leq 3.} 
\end{equation}
\rrred{To start with,} simple calculations show that the function 
\begin{equation}
\label{eq529aabb9}
\rred{\frac{t^2}{2}+F(\lm t)-\mu t\text{ is strictly convex for all }\mu\in \R}
\end{equation}
and hence \rblue{$\frac{t^2}{2}+F(\lm t)=\mu t$} has at most two solutions, with $t=0$ being trivial.  \blue{If $\hat w\times\hat h=0$, then there exists $\mu\not=0$ such that $\hat{w}=\mu \hat{h}$ and in the same way as before,  by (\ref{eq529aabb9}) and (\ref{eq529aabb6}), we either have $h_i=0$ for some $i$ or $h_{i_0}=h_{i_1}$ and thus \rrred{(\ref{ss44}) follows}.}

So \rrred{by (\ref{eq529aabb7}),  (\ref{ss44})} we may assume
	\begin{equation}
	\label{ep5}
	\hat z\times\hat h\ne 0, \qd \hat w\times \hat h\ne 0.
	\end{equation}
	Now the three minors in $\mathrm{Span}\{\Ub_{\KKI}^{0}\}$ are
	\begin{equation}
	\label{ss46}
	M_1 = (\hat h\cdot \hat\alpha)^2-\lm(\hat h\cdot\hat\alpha)(\hat z\cdot\hat\alpha),
	\end{equation}
	\begin{equation}
	\label{ss47}
	M_2 = (\hat h\cdot\hat\alpha)(\hat w\cdot\hat\alpha)-\lm(\hat h\cdot\hat\alpha)(\hat y\cdot\hat\alpha),
	\end{equation}
	and
	\begin{equation}
	\label{ss48}
	M_3 = (\hat z\cdot\hat\alpha)(\hat w\cdot\hat\alpha)-(\hat h\cdot\hat\alpha)(\hat y\cdot\hat\alpha). 
	\end{equation}
	To solve for the Rank-$1$ directions,  from $M_1=0$, we need $\hat h\cdot\hat\alpha=0$ or $\hat h\cdot\hat\alpha=\lm \hat z\cdot\hat\alpha$. When $\hat h\cdot \hat\alpha=0$, it follows from $M_3=0$ that $\hat z\cdot\hat\alpha=0$ or $\hat w\cdot\hat\alpha=0$, and this produces no Rank-$1$ directions due to (\ref{ep5}). When $\hat h\cdot\hat\alpha=\lm \hat z\cdot\hat\alpha$, we have $(\hat h-\lm\hat z)\cdot\hat \alpha=0$ and there is at most one Rank-$1$ direction since $\hat h-\lm\hat z\overset{\rrred{(\ref{ep5})}}{\ne} 0$. Thus we can apply Lemma \ref{L13} (a).
	
	The case when $\dim\lt(\mathrm{Span}\{\Ub_{\KKI}^{0}\}\rt)=3$ can be argued in the same manner as in Step 1 following the above lines. \rrred{We obtain an analogue of (\ref{ss43}) where $\mathrm{Span}\{\Ub_{\KKI}^{0}\}$ is a three-dimensional subspace \bblue{parameterized} by $\vec \alpha\in \R^3$. By exactly the same argument we used to \bblue{establish}
(\ref{eq529aabb7}) and (\ref{ss44}), we have that $\vec z\times \vec h\not=0$ and $\vec w\times \vec h\not=0$. We obtain the same set of minors given by (\ref{ss46}), (\ref{ss47}), (\ref{ss48}). Now \bblue{$M_1=0$ implies} $\vec h \cdot \vec \alpha=0$ or $\vec h \cdot \vec \alpha=\lm \vec z \cdot \vec \alpha$. When $\vec h \cdot \vec \alpha=0$, from 
$M_3=0$ we have $\vec z \cdot \vec \alpha=0$ or $\vec w \cdot \vec \alpha=0$ and so the \rblue{Rank-$1$ directions form two lines}. When $\lt(\vec h-\lm \vec z\rt)\cdot \vec \alpha=0$, since $\vec h-\lm \vec z\not=0$, the Rank-$1$ directions form \bblue{at most} a two-dimensional plane. \rblue{As the Rank-$1$ line given by $\vec h\cdot\vec \alpha=\vec z\cdot\vec \alpha=0$ is contained in the plane $\lt(\vec h - \lm\vec z\rt)\cdot\vec\alpha=0$,} by Lemma \ref{L13} (b) we are done}. This completes the proof of Step 2 and the lemma.
\end{proof}

%
%

%
%

\section{Non-existence of three-dimensional $\Tb_4$}

In this section we prove non-existence of three-dimensional $\Tb_4$ in $\KKI_1$. We denote
\begin{eqnarray}
	\label{eq92.2}
	\blue{\SI_{\KKI}^0}:=
	\lt(\begin{array}{ccccc} h_1 & r_1 & a(r_1) & h_1a(r_1) & \frac{h_1^2}{2}+F(r_1)  \\ 
		h_2 & r_2 & a(r_2) & h_2a(r_2) & \frac{h_2^2}{2}+F(r_2)\\ 
		h_3 & r_3 & a(r_3) & h_3a(r_3) & \frac{h_3^2}{2}+F(r_3)
	\end{array}\rt).
\end{eqnarray}

\begin{a1}
	\label{L1.22}
	Let $\Ub_{\KKI}^{0}$ be defined by (\ref{eqe25}) \green{and $\SI_{\KKI}^0$ be defined by (\ref{eq92.2})}, then 
	\begin{equation*} 
		\mathrm{Rank}(\SI_{\KKI}^0)=p\Longleftrightarrow \mathrm{dim}\lt( \mathrm{Span}\{\Ub_{\KKI}^{0}\} \rt)=p \qd\text{ for }p=2,3.
	\end{equation*}
\end{a1}
\begin{proof}
	Writing out the entries of $\QI(h_i,r_i)$ as the rows of a matrix we have 
	that 
	\begin{equation*}
		\mathrm{dim}\lt( \mathrm{Span}\lt\{\QI(h_i,r_i): i=1,2,3\rt\}\rt)=p
	\end{equation*}
	is equivalent to 
	\begin{equation*}
		\mathrm{Rank}\lt(\begin{array}{cccccc} h_1 &  a(r_1) & h_1a(r_1) & r_1 & h_1 & \frac{h_1^2}{2}+F(r_1)  \\ 
			h_2 &  a(r_2) & h_2a(r_2) & r_2 & h_2 & \frac{h_2^2}{2}+F(r_2)  \\ 
			h_3 &  a(r_3) & h_3a(r_3) & r_3 & h_3 & \frac{h_3^2}{2}+F(r_3)  
		\end{array}\rt)=p. 
	\end{equation*}
	It is immediate that this is equivalent to $\mathrm{Rank}(\SI_{\KKI}^0)=p$ for $p=2,3$. 
\end{proof} 

Our main result of this section is
\begin{a2}
	\label{T3}
	Let $\Ub_{\KKI}^{0}$ be defined by (\ref{eqe25}). If 
	$\mathrm{dim}\lt(\mathrm{Span}\{\Ub_{\KKI}^{0}\}\rt)=3$ then $\KKI$ cannot \blue{contain} a $\Tb_4$.
\end{a2}
The proof is done in several steps. To this end we define
\begin{equation*}
\AI^l_{\KKI}:=\lt(\begin{array}{ccc} h_1 & h_2 & h_3\\
a(r_1) & a(r_2) & a(r_3)\\
h_1 a(r_1) & h_2 a(r_2) & h_3 a(r_3)\end{array}\rt)
\end{equation*}
and
\begin{equation*}
\AI^r_{\KKI}:=\lt(\begin{array}{ccc} r_1 & r_2 & r_3\\
h_1 & h_2 & h_3\\
\frac{h_1^2}{2}+F(r_1) & \frac{h_2^2}{2}+F(r_2) & \frac{h_3^2}{2}+F(r_3)\end{array}\rt).
\end{equation*}
Further we denote
\begin{equation*}
\AI_{\KKI}^0 := \lt(\begin{array}{cc}\AI^l_{\KKI} &\AI^r_{\KKI}\end{array}\rt)\in M^{3\times 6}.
\end{equation*}

\begin{a1}
\label{L12}
Assume $\dim\lt(\mathrm{Span}\{\Ub_{\KKI}^{0}\}\rt)=3$. If $\mathrm{Rank}(\AI^l_{\KKI})=3$ or $\mathrm{Rank}(\AI^r_{\KKI})=3$, then $\Ub_{\KKI}^{0}$ cannot \blue{contain} a $\Tb_4$.
\end{a1}

\begin{proof}
Without loss of generality, assume $\mathrm{Rank}(\AI^l_{\KKI})=3$. The case when $\mathrm{Rank}(\AI^r_{\KKI})=3$ can be dealt with in exactly the same manner. Note that the subspace 
\begin{equation*}
\WI:=\mathrm{Span}\{\QI(h_1,r_1),  \QI(h_2,r_2),\QI(h_3,r_3) \}
\end{equation*}
can be parameterized by the mapping $\PPI:\R^3\rightarrow M^{3\times 2}$ defined by 
\begin{equation}
\label{eq105.4}
\PPI(x):=
\lt(\begin{array}{cc} \AI^l_{\KKI} x & \AI^r_{\KKI} x\end{array}\rt).
\end{equation}
Denote by $[\AI^{*}_{\KKI}]_k$ the $k$-th row of the matrix $\AI^{*}_{\KKI}$. As $\mathrm{Rank}(\AI^l_{\KKI})=3$, the three rows of $\AI^l_{\KKI}$ are linearly independent. Hence we can write
\begin{equation*}
[\AI^r_{\KKI}]_j = \sum_{k=1}^3 \lm_{jk}[\AI^l_{\KKI}]_k
\end{equation*}
for some $\lm_{jk}\in\R$. Denoting the matrix \blue{$\BI\in M^{3\times 3}$ by $[\BI]_{jk}:=\lm_{jk}$}, we have
\begin{equation*}
\BI \AI^l_{\KKI} = \AI^r_{\KKI}
\end{equation*}
and it follows that
\begin{equation}
\label{eqp14}
\AI^r_{\KKI} x = \BI \lt(\AI^l_{\KKI} x\rt)
\end{equation}
\blue{for all $x\in\R^3$}. So letting $y:= \AI^l_{\KKI} x$, it follows that
\begin{equation}\label{eq138}
\PPI(x)=\PPI\lt(\lt(\AI^l_{\KKI}\rt)^{-1}y\rt)\overset{(\ref{eq105.4}),(\ref{eqp14})}{=}\lt(\begin{array}{cc} y & \BI y \end{array}\rt).
\end{equation}

Now $\mathrm{Rank}( \PPI((\AI^l_{\KKI})^{-1} y))=1$ if and only if $y$ is an eigenvector of the matrix $\BI$. So there 
are three possibilities to consider:  $\BI$ either has one, two or three distinct eigenvalues. If $\BI$ has three distinct eigenvalues, since the dimension of the eigenspace is bounded above by the multiplicity of 
the corresponding eigenvalue,  $\BI$ has three \green{linearly independent} eigenvectors and thus
$\Lambda_R\cap \WI$ consists of three distinct lines. If $\BI$ has
two distinct eigenvalues, \blue{the dimensions of the eigenspaces are either two and one or one and one.} Therefore $\Lambda_R\cap \WI$ either 
consists of two distinct lines, or a two-dimensional plane and a line. So in the above cases, it follows from Lemma \ref{L13} that $\Ub_{\KKI}^{0}$ cannot \blue{contain} a $\Tb_4$.

Finally suppose $\BI$ has just one eigenvalue. If 
the dimension of the eigenspace is less than three then the situation reduces to the ones already discussed and the conclusion of the lemma follows. So suppose the dimension of the eigenspace is three, then every vector $y$ is an eigenvector of $\BI$ \blue{and from (\ref{eq138}) we immediately have $\mathrm{Rank}(\PPI(x))=1$ for all $x\in\R^3$. As $\dim\lt(\WI\rt)=3$, it is clear that $\PPI:\R^3\rightarrow\WI$ is a linear isomorphism and hence $\WI\subset\Lambda_{R}$. In particular, $\QI(h_i,r_i)-\QI(0,0)=\QI(h_i,r_i)\in \Lambda_{R}$ and thus $\Ub_{\KKI}^{0}$ contains Rank-$1$ connections. By definition, $\Ub_{\KKI}^{0}$ is not a $\Tb_4$.}
\end{proof}

\begin{a1}
	\label{L15}
	Assume $\dim\lt(\mathrm{Span}\{\Ub_{\KKI}^{0}\}\rt)=3$. If $\mathrm{Rank}(\AI^{*}_{\KKI})=2$ for $*=l,r$ and $\mathrm{Rank}(\AI_{\KKI}^0)=3$, then $\Ub_{\KKI}^{0}$ cannot \blue{contain} a  $\Tb_4$.
\end{a1}

\begin{proof}
	\green{By Lemma \ref{L1.5}} we can assume that $\vec h\times\vec r\ne 0$ and $\vec h\times\vec z\ne 0$ (recall that $\vec z$ is defined in (\ref{eqp2})). As $\mathrm{Rank}(\AI^l_{\KKI})=\mathrm{Rank}(\AI^r_{\KKI})=2$, there exist $\lm_1, \lm_2, \mu_1, \mu_2$ such that
\begin{equation}
\label{eqp21}
h_i a(r_i) = \lm_1 h_i + \lm_2 a(r_i)
\end{equation}
and
\begin{equation}
\label{eqp22}
\frac{h_i^2}{2}+F(r_i) = \mu_1 r_i + \mu_2 h_i. 
\end{equation}
Therefore we have
\begin{equation*}
\mathrm{Span}\{\Ub_{\KKI}^{0}\}=\lt\{\lt(\begin{array}{cc} \vec h\cdot\vec\alpha & \vec r\cdot\vec\alpha \\ \vec z\cdot\vec\alpha & \vec h\cdot\vec\alpha \\ \lm_1(\vec h\cdot\vec\alpha) + \lm_2 (\vec z\cdot\vec\alpha) & \mu_1 (\vec r\cdot\vec\alpha) + \mu_2 (\vec h\cdot\vec\alpha)\end{array}\rt): \vec\alpha\in\R^3\rt\}.
\end{equation*}
Since $\mathrm{Rank}(\AI_{\KKI}^0)=3$, \green{we must have 
\begin{equation}
\label{eqbbvv6}
(\lm_1,\lm_2)\ne (\mu_1,\mu_2), 
\end{equation}
as} otherwise the third row of $\AI_{\KKI}^0$ would be a linear combination of the first two rows of $\AI_{\KKI}^0$, which contradicts our assumption. 
	
	Now we calculate the three minors in $\mathrm{Span}\{\Ub_{\KKI}^{0}\}$ and get
	\begin{equation*}
	M_1 = (\vec h\cdot\vec\alpha)^2 -(\vec r\cdot\vec\alpha)(\vec z\cdot\vec\alpha),
	\end{equation*}
	\begin{equation*}
	\begin{split}
	M_2 &= (\vec h\cdot\vec\alpha)\lt(\mu_1 (\vec r\cdot\vec\alpha) + \mu_2 (\vec h\cdot\vec\alpha)\rt) - (\vec r\cdot\vec\alpha)\lt(\lm_1(\vec h\cdot\vec\alpha) + \lm_2 (\vec z\cdot\vec\alpha)\rt)\\
	&= \mu_1 (\vec h\cdot\vec\alpha)(\vec r\cdot\vec\alpha) + \mu_2 (\vec h\cdot\vec\alpha)^2 - \lm_1 (\vec h\cdot\vec\alpha)(\vec r\cdot\vec\alpha) -\lm_2 (\vec r\cdot\vec\alpha)(\vec z\cdot\vec\alpha)\\
	&= (\mu_1-\lm_1)(\vec h\cdot\vec\alpha)(\vec r\cdot\vec\alpha) + (\mu_2-\lm_2)(\vec h\cdot\vec\alpha)^2 + \lm_2 \lt((\vec h\cdot\vec\alpha)^2-(\vec r\cdot\vec\alpha)(\vec z\cdot\vec\alpha)\rt),
	\end{split}
	\end{equation*}
	and
	\begin{equation*}
	\begin{split}
	M_3 &= (\vec z\cdot\vec\alpha)\lt(\mu_1 (\vec r\cdot\vec\alpha) + \mu_2 (\vec h\cdot\vec\alpha)\rt) - (\vec h\cdot\vec\alpha)\lt(\lm_1(\vec h\cdot\vec\alpha) + \lm_2 (\vec z\cdot\vec\alpha)\rt)\\
	&= \mu_1 (\vec z\cdot\vec\alpha)(\vec r\cdot\vec\alpha) + \mu_2 (\vec z\cdot\vec\alpha)(\vec h\cdot\vec\alpha) - \lm_1 (\vec h\cdot\vec\alpha)^2 -\lm_2 (\vec h\cdot\vec\alpha)(\vec z\cdot\vec\alpha)\\
	&= (\mu_1-\lm_1)(\vec h\cdot\vec\alpha)^2 + (\mu_2-\lm_2)(\vec h\cdot\vec\alpha)(\vec z\cdot\vec\alpha) + \mu_1\lt((\vec r\cdot\vec\alpha)(\vec z\cdot\vec\alpha)-(\vec h\cdot\vec\alpha)^2\rt).
	\end{split}
	\end{equation*}
	
	When $\vec h\cdot\vec\alpha=0$, the Rank-$1$ directions must satisfy $M_1=0$ and so we need $(\vec r\cdot\vec\alpha)(\vec z\cdot\vec\alpha)=0$. Recall that $\vec h$ and $\vec r$, $\vec h$ and $\vec z$ are both linearly independent. When $\vec h\cdot\vec\alpha=\vec r\cdot \vec\alpha=0$, we get one Rank-$1$ direction. There is another Rank-$1$ direction when $\vec h\cdot\vec\alpha=\vec z\cdot \vec\alpha=0$. When $\vec h\cdot\vec\alpha\ne 0$, $M_1=M_2=M_3=0$ is equivalent to (note that $M_1$ is part of the expressions in the last lines of $M_2$ and $M_3$)
	\begin{equation*}
	(\vec h\cdot\vec\alpha)^2 -(\vec r\cdot\vec\alpha)(\vec z\cdot\vec\alpha)=0,
	\end{equation*}
	\begin{equation}
	\label{eqp19}
	(\mu_1-\lm_1)(\vec r\cdot\vec\alpha) + (\mu_2-\lm_2)(\vec h\cdot\vec\alpha)=0
	\end{equation}
	and
	\begin{equation}
	\label{eqp20}
	(\mu_1-\lm_1)(\vec h\cdot\vec\alpha) + (\mu_2-\lm_2)(\vec z\cdot\vec\alpha)=0.
	\end{equation}
	\blue{Thus, rewriting (\ref{eqp19}) and (\ref{eqp20}), the Rank-$1$ directions must satisfy
	\begin{equation}
	\label{eqp23}
	\lt((\mu_1-\lm_1)\vec r+(\mu_2-\lm_2)\vec h\rt) \cdot \vec\alpha=0 \text{ and }  \lt((\mu_1-\lm_1)\vec h+(\mu_2-\lm_2)\vec z\rt) \cdot \vec\alpha=0.
	\end{equation}
	Since $\dim\lt(\mathrm{Span}\{\Ub_{\KKI}^0\}\rt)=3$, we know from Lemma \ref{L1.22} that $\mathrm{Rank}(\SI_{\KKI}^0)=3$. Equations (\ref{eqp21}) and (\ref{eqp22}) imply that \rblue{all the column vectors of $\SI_{\KKI}^0$ are in $\mathrm{Span}\{\vec h,\vec r,\vec z\}$}, and hence $\vec h,\vec r,\vec z$ must be linearly independent. Because of (\ref{eqbbvv6}), we must have $\lm_1\ne\mu_1$ or $\lm_2\ne\mu_2$, and it follows immediately that  $(\mu_1-\lm_1)\vec r+(\mu_2-\lm_2)\vec h$ and $(\mu_1-\lm_1)\vec h+(\mu_2-\lm_2)\vec z$ are linearly independent. Hence (\ref{eqp23}) gives only one possible Rank-$1$ line in the case when $\vec h\cdot\vec\alpha\ne 0$. Combining this with the case $\vec h\cdot\vec\alpha=0$ we see that there are at most three distinct Rank-$1$ directions in $\mathrm{Span}\{\Ub_{\KKI}^0\}$. An application of Lemma \ref{L13} (b) completes the proof.}
\end{proof}

\subsection{The case $\mathrm{Rank}(\AI_{\KKI}^0)=2$}\label{subsec:rank2}
It only remains to consider the case when $\mathrm{Rank}(\AI_{\KKI}^0)=2$. We need a key lemma concerning the function $F$. \zred{We state the result in more general form for later application in Proposition \ref{TCXE}. 
\begin{a1}
\label{l18}
Suppose $I\bbblue{\subset\R}$ is an open interval containing $0$. Let $\ti{a}\in C^2(I)$ \bbblue{be such that $\ti a'>0$ and $\ti a(0)=0$}. Suppose $\wt{F}$ is a primitive of $\ti{a}$ \bbblue{with $\wt{F}(0)=0$}. If 
$\ti{a}$ is strictly convex \bbblue{in $I$}  then
\begin{equation}
\label{eqp28}
2\wt{F}(r)-r\ti{a}(r)>0 \text{ for } r<0 \qd\text{ and }\qd 2\wt{F}(r)-r\ti{a}(r)<0 \text{ for } r>0.
\end{equation}
And if $\ti{a}$ is strictly concave \bbblue{in $I$} then
\begin{equation*}
2\wt{F}(r)-r\ti{a}(r)<0 \text{ for } r<0 \qd\text{ and }\qd 2\wt{F}(r)-r\ti{a}(r)>0 \text{ for } r>0.
\end{equation*}
\end{a1}}
\begin{proof}
\zred{We argue only in the case where  $\ti{a}$ is strictly convex; the case where  $\ti{a}$ is strictly concave follows in the same way. Letting $g(r):=2\wt{F}(r)-r\ti{a}(r)$ and using $\wt{F}'(r)=\ti{a}(r)$, we have $g'(r)=\ti{a}(r)-r\ti{a}'(r)$ and $g''(r)=-r\ti{a}''(r)$. Since $\ti{a}$ is strictly convex, we know $g''>0$ for $r<0$ and $g''<0$ for $r>0$. Further, as $\wt F(0)=\ti a(0)=0$, we know $g(0)=0$ and $g'(0)=0$. Combining this with the sign for $g''$, we know that $g'(r)<0$ for $r<0$ and $g'(r)<0$ for $r>0$. It follows that $g(r)>0$ for $r<0$ and $g(r)<0$ for $r>0$ and this translates to exactly (\ref{eqp28}).}
\end{proof}

As before, we may assume that $\vec h\times\vec r\ne 0$ and $\vec h\times \vec z\ne 0$ where $\vec z=(a(r_1),a(r_2),a(r_3))$, \blue{and thus the first two rows of $\AI_{\KKI}^0$ are linearly independent}. So there exist $\lm_1$ and $\lm_2$ such that
\begin{equation}
\label{eqp26}
h_i a(r_i) = \lm_1 h_i + \lm_2 a(r_i)
\end{equation}
and
\begin{equation}
\label{eqp27}
\frac{h_i^2}{2}+F(r_i) = \lm_1 r_i + \lm_2 h_i. 
\end{equation}
We define the set
\begin{equation*}
\ti\Ub_{\KKI}^{0}:=\lt\{\lt(\begin{array}{cc}0 & 0\\ 0& 0\end{array}\rt), \lt(\begin{array}{cc}h_1& r_1\\ a(r_1)& h_1\end{array}\rt),\lt(\begin{array}{cc}h_2 & r_2\\ a(r_2)& h_2\end{array}\rt),\lt(\begin{array}{cc}h_3 & r_3\\ a(r_3)& h_3\end{array}\rt)\rt\}.
\end{equation*}

\begin{a1}
	\label{L100}
	\blue{If $(h_i,r_i)$ satisfies the system (\ref{eqp26})-(\ref{eqp27}) for $i=1,2,3$ and $\Ub_{\KKI}^{0}$ with the given ordering forms a $\Tb_4$ with $\dim\lt(\mathrm{Span}\{\Ub_{\KKI}^{0}\}\rt)=3$, then $\ti\Ub_{\KKI}^{0}$ also forms a $\Tb_4$ with the given ordering.}
\end{a1}
\begin{proof}
	\blue{We define the linear mapping $L: \mathrm{Span}\{\Ub_{\KKI}^{0}\}\rightarrow\mathrm{Span}\{\ti\Ub_{\KKI}^{0}\}$ by $$L\lt(\QI(h_i,r_i)\rt)=\lt(\begin{array}{cc}h_i& r_i\\ a(r_i)& h_i\end{array}\rt) \qd\text{ for }i=1,2,3.$$
	\rblue{Noting (\ref{eqp26})-(\ref{eqp27}), it is clear that} $L$ satisfies (\ref{eqa30}). Since $\dim\lt(\mathrm{Span}\{\Ub_{\KKI}^{0}\}\rt)=3$, we know  $\mathrm{Rank}(\SI_{\KKI}^0)=3$ from Lemma \ref{L1.22}, and thus $\vec h, \vec r, \vec z$ are linearly independent because of (\ref{eqp26})-(\ref{eqp27}). Thus $\dim\lt(\mathrm{Span}\{\ti\Ub_{\KKI}^{0}\}\rt)=3$ and therefore the mapping $L$ is a linear isomorphism. Now Lemma \ref{LL0} applies to finish the proof.}
\end{proof}

\begin{a1}
	\label{L17}
	Assume $\dim\lt(\mathrm{Span}\{\Ub_{\KKI}^{0}\}\rt)=3$. If $(h_i,r_i)$ satisfies the system (\ref{eqp26})-(\ref{eqp27}) with $\lm_1=0$ or $\lm_2=0$ for $i=1,2,3$, then $\Ub_{\KKI}^{0}$ cannot \blue{contain} a $\Tb_4$.
\end{a1}

\begin{proof}
	By Lemma \ref{L100}, it suffices to show the set $\ti\Ub_{\KKI}^{0}$ cannot \blue{contain} a $\Tb_4$. Our main tool is Proposition \ref{rmk2}. Recall that from (\ref{eq8}) the function $F$ is strictly convex with $F(0)=F'(0)=0$, and thus $F\geq 0$ for all $r$ and $F=0$ only at $r=0$. If $\lm_1=\lm_2=0$, from (\ref{eqp27}) we must have $h_i=0$ and $r_i=0$ for all $i=1,2,3$, in which case \blue{$\ca{\ti\Ub_{\KKI}^{0}}<4$}. So we only have to consider the case when $\lm_1\ne 0$ or $\lm_2\ne 0$. \nl
	
\em Step 1. \rm We first consider the case when \red{$\lm_1=0, \blue{\lm_2\ne 0}$.}
So (\ref{eqp26})-(\ref{eqp27}) become
\begin{equation}\label{ep10}
h_ia(r_i)=\lm_2a(r_i),
\end{equation}
\begin{equation}\label{ep11}
\frac{h_i^2}{2}+F(r_i)=\lm_2 h_i.
\end{equation}
From (\ref{ep10}), we have 
\begin{equation}
\label{cgeq2}	
\red{h_i=\lm_2\text{ for any }i\text{ for which }a(r_i)\not=0.}
\end{equation}
\rrred{Let $\Pi_r:=\lt\{i\bblue{\in\{1,2,3\}}:r_i\not=0\rt\}$. So $\ca{\Pi_r}\in \lt\{0,1,2,3\rt\}$. We consider each case in turn.} 

\rrred{Case 1: $\ca{\Pi_r}=3$. So} by (\ref{eq7}) we have $a(r_i)\ne 0$ for $i=1,2,3$. Thus $h_i=\lm_2$ for all $i$. Then from (\ref{ep11}), $r_i$ solves
\begin{equation}
\label{eqp29}
F(r)=\frac{\lm_2^2}{2}
\end{equation}
for all $i$. But as $F$ is strictly convex, this equation has at most two distinct roots and hence $\ca{\ti\Ub_{\KKI}^{0}}\leq 3$. 

\rrred{Case 2: $\ca{\Pi_r}\leq 1$. So} there exist $i\ne j$ such that $r_i=r_j=0$, then from (\ref{ep11}), $h_i$ and $h_j$ both solve the equation $\frac{h^2}{2}-\lm_2 h=0$ which has the solutions $0$ and $2\lm_2$. If $h_i=0$ or $h_j=0$, then $(h_i,r_i)=(0,0)$ or $(h_j,r_j)=(0,0)$. Otherwise, we have $(h_i,r_i)=(h_j,r_j)=(2\lm_2,0)$. In both cases, $\ca{\ti\Ub_{\KKI}^{0}}\leq 3$. 

\rrred{Case 3: $\ca{\Pi_r}=2$. So} exactly one of the $r_i$'s equals zero, without loss of generality, \red{assume $r_1=0$.} From (\ref{ep11}), we know 
\begin{equation}
\label{cgeq4}
\red{h_1=2\lm_2.}
\end{equation}
(\blue{Otherwise $(h_1,r_1)=(0,0)$ and $\ca{\ti\Ub_{\KKI}^{0}}\leq 3$.}) Also we have 
\begin{equation}
\label{cgeq5}
\red{h_2=h_3\overset{(\ref{cgeq2})}{=}\lm_2.}
\end{equation}
Then $r_2$ and $r_3$ are solutions of (\ref{eqp29}), which has at most two distinct solutions. If (\ref{eqp29}) fails to have two distinct solutions, then we are done. So assume that (\ref{eqp29}) has two distinct solutions, \blue{then $r_2$ and $r_3$ must take these two distinct solutions in order for $\ca{\ti\Ub_{\KKI}^{0}}=4$. Because of (\ref{eq8}), the two distinct solutions of (\ref{eqp29}) must have opposite signs. Without loss of generality}, assume $r_2<0<r_3$. From Lemma \ref{l18} we have
	\begin{equation}
	\label{eqp31}
	\lm_2^2-r_3a(r_3)\overset{(\ref{eqp29})}{=} 2F(r_3)-r_3a(r_3)<0.
	\end{equation}
	Now our set $\ti\Ub_{\KKI}^{0}$ becomes
	\begin{equation*}
	\ti\Ub_{\KKI}^{0}\red{\overset{(\ref{cgeq4}),(\ref{cgeq5})}{=}}\lt\{\lt(\begin{array}{cc}0 & 0\\ 0& 0\end{array}\rt),\lt(\begin{array}{cc}2\lm_2 & 0\\ 0& 2\lm_2\end{array}\rt),\lt(\begin{array}{cc}\lm_2 & r_2\\ a(r_2)& \lm_2\end{array}\rt),\lt(\begin{array}{cc}\lm_2 & r_3\\ a(r_3)& \lm_2\end{array}\rt)\rt\}.
	\end{equation*}
	We call the above matrices $T_0, T_1, T_2, T_3$. Now we observe that
	\begin{equation*}
	\blue{\det(T_0-T_3)=\det(T_1-T_3)=\lm_2^2-r_3a(r_3)\red{\overset{(\ref{eqp31})}{<}}0,}
	\end{equation*}
	and
	\begin{equation*}
	\blue{\det(T_2-T_3)=-(r_3-r_2)(a(r_3)-a(r_2))<0,}
	\end{equation*}
	\blue{where the last inequality holds because the function $a$ is strictly increasing.} Since $\det(T_i-\blue{T_3})\blue{<}0$ for all $i\ne 3$, it follows from Proposition \ref{rmk2} that $\ti\Ub_{\KKI}^{0}$ cannot \blue{contain} a $\Tb_4$. This completes \rrred{the proof of Case 3 and the proof of Step 1.} \nl
	
	\em Step 2. \rm Next we consider the case when $\lm_2=0$ and \blue{$\lm_1\ne 0$}.  Now (\ref{eqp26})-(\ref{eqp27}) become
	\begin{equation}\label{ep12}
	h_ia(r_i)=\lm_1h_i,
	\end{equation}
	\begin{equation}\label{ep13}
	\frac{h_i^2}{2}+F(r_i)=\lm_1 r_i.
	\end{equation}
	From (\ref{ep12}) we know $a(r_i)=\lm_1$ unless $h_i=0$. \rrred{Similarly to how we argued in Step 1 we let $\Pi_h:=\lt\{i\bblue{\in\{1,2,3\}}:h_i\not=0\rt\}$. So $\ca{\Pi_h}\in \lt\{0,1,2,3\rt\}$. Again we consider each case in turn.} 

\rrred{Case 1: $\ca{\Pi_h}=3$. So }$h_i\ne 0$ for all $i$ and we have $a(r_i)=\lm_1$ for all $i$. As $a$ is strictly monotonic, this implies that all $r_i$'s are equal, and hence from (\ref{ep13}) all $h_i^2$ equals the same constant. It is a simple argument to see that $\ca{\ti\Ub_{\KKI}^{0}}<4$ in this case. 

\rrred{Case 2: $\ca{\Pi_h}\leq 1$. So} $h_i=h_j=0$ for some $i\ne j\in\{1,2,3\}$, and it follows from (\ref{ep13}) that $r_i$ and $r_j$ both solve $F(r)=\lm_1 r$, which has at most one non-trivial solution. \rrred{As in Case 2 of Step 1} it is easy to see that $\ca{\ti\Ub_{\KKI}^{0}}<4$ in this case. 
	
\rrred{Case 3: $\ca{\Pi_h}=2$. So} exactly one of the $h_i$'s vanishes. Without loss of generality, assume $h_1=0$. It follows from \rrred{(\ref{ep13})} that $r_1$ must be the non-trivial solution of $F(r)=\lm_1 r$ in order for $\ca{\ti\Ub_{\KKI}^{0}}=4$. As $h_2\ne 0$ and $h_3\ne 0$, \rrred{from (\ref{ep12})} we have $a(r_2)=a(r_3)=\lm_1$ and hence $r_2=r_3=:\sigma$. From (\ref{ep13}), $h_2$ and $h_3$ solve $h^2=2\lm_1 \sigma-2F(\sigma)$. So \bblue{this equation must have} two distinct solutions (if not then $\ca{\ti\Ub_{\KKI}^{0}}<4$), and denote them by $-\beta, \beta$. Without loss of generality, we let $h_2=-\beta$ and $h_3=\beta$. Now we have
	\begin{equation*}
	\ti\Ub_{\KKI}^{0}=\lt\{\lt(\begin{array}{cc}0 & 0\\ 0& 0\end{array}\rt),\lt(\begin{array}{cc}0 & r_1\\ a(r_1)& 0\end{array}\rt),\lt(\begin{array}{cc}-\beta & \sigma\\ a(\sigma)& -\beta\end{array}\rt),\lt(\begin{array}{cc}\beta & \sigma\\ a(\sigma)& \beta\end{array}\rt)\rt\}.
	\end{equation*}
	As in Step 1, we label the matrices in $\ti\Ub_{\KKI}^{0}$ by $T_0,T_1,T_2,T_3$ and calculate
	\begin{equation*}
	\det(T_1-T_0)=-r_1a(r_1)<0,\qd \det(T_2-T_0)=\beta^2-\sigma a(\sigma), 
	\end{equation*}
	\begin{equation*}
	\det(T_3-T_0)=\beta^2-\sigma a(\sigma),\qd \det(T_2-T_1)=\beta^2-(\sigma-r_1)(a(\sigma)-a(r_1)),\end{equation*}
	\begin{equation*}
	\det(T_3-T_1)=\beta^2-(\sigma-r_1)(a(\sigma)-a(r_1)),\qd\det(T_3-T_2)=4\beta^2>0. 
	\end{equation*}
	We denote $d_1:=\beta^2-\sigma a(\sigma)$ and $d_2:=\beta^2-(\sigma-r_1)(a(\sigma)-a(r_1))$.  If $d_1<0$, then $\det(T_i-T_0)<0$ for all $i\ne 0$. \red{If $d_2<0$}, then $\det(T_i-T_1)<0$ for all $i\ne 1$.  If $d_1>0$ and $d_2>0$, then $\det(T_i-T_3)>0$ for all $i\ne 3$. In conclusion, we can always find some $T_i$ such that $\{\det(T_j-T_i)\}$ does not change sign. Again by Proposition \ref{rmk2}, $\ti\Ub_{\KKI}^{0}$ cannot \blue{contain} a $\Tb_4$ and this completes Step 2. 
	\end{proof}

It remains to consider the case when $(h_i,r_i)$ satisfies (\ref{eqp26})-(\ref{eqp27}) for $i=1,2,3$ with $\lm_1\ne 0$ and $\lm_2\ne 0$. We collect some elementary facts about the system (\ref{eqp26})-(\ref{eqp27}). \blue{First note that if $a(r_i)=\lm_1$ for some $i$, then equation (\ref{eqp26}) would imply $\lm_2a(r_i)=0$. This would yield $\lm_1=a(r_i)=0$ which is a contradiction. So we must have}
\begin{equation}
\label{ss63}
a(r_i)\ne\lm_1\qd\text{ for }i=1,2,3.
\end{equation}
	
	\begin{a1}\label{l20}
		The system (\ref{eqp26})-(\ref{eqp27}) has at most two distinct solutions satisfying $a(r)<\lm_1$.
	\end{a1}

	\begin{proof}
		Let $(h,r)$ be a solution to the system (\ref{eqp26})-(\ref{eqp27}). We can solve for $h$ from (\ref{eqp26}) and get
		\begin{equation}
		\label{eqp32}
		h=\frac{\lm_2 a(r)}{a(r)-\lm_1}.
		\end{equation}
		Plugging this into (\ref{eqp27}) we obtain that $r$ solves
		\begin{equation*}
		\frac{\lm_2^2 a(r)^2}{2(a(r)-\lm_1)^2} + F(r)=\lm_1 r+\frac{\lm_2^2 a(r)}{a(r)-\lm_1}.
		\end{equation*}
		Simplifying the above equation, we obtain
		\begin{equation*}
		F(r)-\lm_1 r-\frac{\lm_2^2}{2} = -\frac{\lm_1^2\lm_2^2}{2(a(r)-\lm_1)^2}.
		\end{equation*}
		Let us denote
		\begin{equation*}
		p(r):=F(r)-\lm_1 r-\frac{\lm_2^2}{2}
		\end{equation*}
		and 
		\begin{equation*}
		q(r):=-\frac{\lm_1^2\lm_2^2}{2(a(r)-\lm_1)^2}.
		\end{equation*}
		Direct calculations \blue{using $F'=a$} show that
		\begin{equation*}
		p'(r)=a(r)-\lm_1, \qd p''(r)=a'(r),
		\end{equation*}
		and
		\begin{equation*}
		q'(r)=\frac{\lm_1^2\lm_2^2a'(r)}{(a(r)-\lm_1)^3}, \qd q''(r)=\lm_1^2\lm_2^2\frac{a''(r)\lt(a(r)-\lm_1\rt)-3a'(r)^2}{(a(r)-\lm_1)^{\blue{4}}}.
		\end{equation*}
		Since $a'(r)>0$, the function $p(r)$ is always strictly convex. For $a(r)<\lm_1$, we have $a''(r)(a(r)-\lm_1)<0$ and thus $q''(r)<0$ for $a(r)<\lm_1$. So the functions $p$ and $q$ can intersect at most twice for $a(r)<\lm_1$ as $q$ is strictly concave here. This completes the proof of the lemma.
	\end{proof}

\begin{a1}\label{l21}
	\blue{Let $(h,r)$ be a non-trivial solution of the system (\ref{eqp26})-(\ref{eqp27}) with $h^2-ra(r)\ne 0$.} If $\lm_1>0$ and $a(r)<\lm_1$, then $h^2-ra(r)>0$; on the other hand, if $\lm_1<0$ and $a(r)>\lm_1$, then $h^2-ra(r)<0$.
\end{a1}

\begin{proof} \rrred{First note that \bblue{$r\ne 0$, as otherwise it} follows from (\ref{eqp26}) and (\ref{eq7}) that $h=0$ and hence $(h,r)$ is a trivial solution \bblue{of (\ref{eqp26})-(\ref{eqp27})}.} We start with $\lm_1>0$. \red{Assume first}
\begin{equation}
\label{cgeq9}	
\red{0<a(r)<\lm_1.}
\end{equation}
It follows from (\ref{eqp32}) that \blue{$\lm_2h<0$}. Solving for $\lm_1,\lm_2$ from (\ref{eqp26}) and (\ref{eqp27}) we obtain
	\begin{equation}
	\label{cgeq12}
	\lm_1 = \frac{a(r)\lt(\frac{h^2}{2}-F(r)\rt)}{h^2-ra(r)}, \qd \lm_2 = \frac{h\lt(\frac{h^2}{2}+F(r)-ra(r)\rt)}{h^2-ra(r)}.
	\end{equation}
	Since $\lm_1>0$ and $a(r)>0$, we know \blue{from the expression for $\lm_1$ that}
	\begin{equation}
	\label{eqp40}
	\lt(\frac{h^2}{2}-F(r)\rt)\lt(h^2-ra(r)\rt)>0.
	\end{equation}
	On the other hand, since \blue{$\lm_2h<0$}, it follows \blue{from the expression for $\lm_2$} that
	\begin{equation}
	\label{eqp41}
	\lt(\frac{h^2}{2}+F(r)-ra(r)\rt)\lt(h^2-ra(r)\rt)<0.
	\end{equation}
	Combining (\ref{eqp40}) with (\ref{eqp41}) gives
	\begin{equation}\label{ep16}
	\lt(\frac{h^2}{2}-F(r)\rt)\lt(\frac{h^2}{2}+F(r)-ra(r)\rt)<0.
	\end{equation}
	\blue{Note that from (\ref{eq7}) we have}
	\begin{equation}\label{ep15}
	\blue{a(r)>0\Longleftrightarrow r>0 \qd\text{and}\qd a(r)<0\Longleftrightarrow r<0.}
	\end{equation}
	Using \blue{(\ref{cgeq9}), (\ref{ep15})} and Lemma \ref{l18} we have
	\begin{equation*}
	\frac{h^2}{2}-F(r) > \frac{h^2}{2}+F(r)-ra(r).
	\end{equation*}
	It follows from this and (\ref{ep16}) that
	\begin{equation*}
	\frac{h^2}{2}-F(r)>0 \qd\text{ and }\qd \frac{h^2}{2}+F(r)-ra(r)<0.
	\end{equation*}
	This together with (\ref{eqp40}) or (\ref{eqp41}) yields $h^2-ra(r)>0$. 
	
	If $a(r)<0$, then from (\ref{eqp32}) we have \blue{$\lm_2h>0$}. Thus 
	\red{from (\ref{cgeq12}), }(\ref{eqp40})-(\ref{eqp41}) become 
	\begin{equation}
	\label{cgeq17}
	\lt(\frac{h^2}{2}-F(r)\rt)\lt(h^2-ra(r)\rt)<0
	\end{equation}
	and 
	\begin{equation*}
	\lt(\frac{h^2}{2}+F(r)-ra(r)\rt)\lt(h^2-ra(r)\rt)>0,
	\end{equation*}
	\rrred{which gives (\ref{ep16}) as before. So} using \blue{(\ref{ep15})} and Lemma \ref{l18} we have
	\begin{equation*}
	\frac{h^2}{2}-F(r) < \frac{h^2}{2}+F(r)-ra(r),
	\end{equation*}
	and hence we must have
	\begin{equation*}
	\frac{h^2}{2}-F(r)<0 \qd\text{ and }\qd \frac{h^2}{2}+F(r)-ra(r)>0.
	\end{equation*}
	It follows \red{from (\ref{cgeq17})} that $h^2-ra(r)>0$. This completes the proof of the first half of the lemma.
	
	\blue{Next we consider $\lm_1<0$ and repeat the above lines.} If $\lm_1<a(r)<0$, then \blue{$\lm_2h<0$} from (\ref{eqp32}) and thus (\ref{eqp40})-(\ref{eqp41}) become
	\begin{equation*}
	\lt(\frac{h^2}{2}-F(r)\rt)\lt(h^2-ra(r)\rt)>0
	\end{equation*}
	and 
	\begin{equation*}
	\lt(\frac{h^2}{2}+F(r)-ra(r)\rt)\lt(h^2-ra(r)\rt)<0.
	\end{equation*}
	Using Lemma \ref{l18} we have
	\begin{equation*}
	\frac{h^2}{2}-F(r) < \frac{h^2}{2}+F(r)-ra(r),
	\end{equation*}
	and hence we must have
	\begin{equation*}
	\frac{h^2}{2}-F(r)<0 \qd\text{ and }\qd \frac{h^2}{2}+F(r)-ra(r)>0.
	\end{equation*}
	It follows that $h^2-ra(r)<0$. 
	
	If $a(r)>0$, then \blue{$\lm_2h>0$} from (\ref{eqp32}) and thus 
	\red{from (\ref{cgeq12})}, (\ref{eqp40})-(\ref{eqp41}) become
	\begin{equation*}
	\lt(\frac{h^2}{2}-F(r)\rt)\lt(h^2-ra(r)\rt)<0
	\end{equation*}
	and 
	\begin{equation*}
	\lt(\frac{h^2}{2}+F(r)-ra(r)\rt)\lt(h^2-ra(r)\rt)>0.
	\end{equation*}
	Using Lemma \ref{l18} we have
	\begin{equation*}
	\frac{h^2}{2}-F(r) > \frac{h^2}{2}+F(r)-ra(r),
	\end{equation*}
	and hence we must have
	\begin{equation*}
	\frac{h^2}{2}-F(r)>0 \qd\text{ and }\qd \frac{h^2}{2}+F(r)-ra(r)<0.
	\end{equation*}
	Thus $h^2-ra(r)<0$. This completes the proof of the lemma.
\end{proof}
	
\begin{a1}\label{l22}
	Let $\lm_1>0$. If $(h_1,r_1)$ and $(h_2,r_2)$ are two non-trivial solutions of the system (\ref{eqp26})-(\ref{eqp27}) with $\lm_1<a(r_1)<a(r_2)$, then $h_1^2-r_1a(r_1)>h_2^2-r_2a(r_2)$.
\end{a1}	

\begin{proof}
	Using (\ref{eqp32}), we have, for $i=1,2$,
	\begin{equation*}
	h_i^2-r_ia(r_i) = \frac{\lm_2^2 a(r_i)^2}{(a(r_i)-\lm_1)^2} - r_ia(r_i). 
	\end{equation*}
	Let us define 
	\begin{equation*}
	l(r):=\frac{\lm_2^2 a(r)^2}{(a(r)-\lm_1)^2} - ra(r).
	\end{equation*}
	When $\lm_1>0$, it is clear that $\frac{a(r)}{a(r)-\lm_1}=1+\frac{\lm_1}{a(r)-\lm_1}$ is decreasing for $a(r)>\lm_1$ and $ra(r)$ is increasing, and thus $l(r)$ is a decreasing function for $a(r)>\lm_1\blue{>0}$. 
\end{proof}	

\blue{To finish the proof in the case when $\mathrm{Rank}(\AI_{\KKI}^0)=2$, we need some preparation. Recall that we fix the set $\KKI\subset\KKI_1$, where $\KKI$ given in (\ref{ep7}) consists of four points parameterized by $(u_i,v_i)$ for $i=0,1,2,3$. Now for $k=0,1,2,3$, we extend the notations in (\ref{eqp1}) by defining 
\begin{equation}
\label{cgeq21}
h_i^k:=u_i-u_k,\qd r_i^k:=v_i-v_k,
\end{equation}
and similar to (\ref{eqe25}) we define the set $\Ub_{\KKI}^k$ associated to the set $\KKI$ with respect to the point $P(u_k,v_k)$ by
\begin{equation}
\label{ep17}
\Ub_{\KKI}^{k}:=\lt\{ \QI_{\blue{v_k}}(h_0^k,r_0^k) , \QI_{\blue{v_k}}(h_1^k,r_1^k), \QI_{\blue{v_k}}(h_2^k,r_2^k),\QI_{\blue{v_k}}(h_3^k,r_3^k)\rt\}.
\end{equation}
Note that when $k=0$, the set $\Ub_{\KKI}^{k}$ agrees with the set $\Ub_{\KKI}^{0}$ defined in (\ref{eqe25}). A crucial observation is that, for $k\in\{1,2,3\}$, we could have switched the labeling of $k$ and $0$ in the set $\KKI$ and thus all the results proved so far also apply to the set $\Ub_{\KKI}^{k}$. Hence it only remains to show}
	
\begin{a1}
	\label{L23}
	\blue{Let $\KKI\subset\KKI_1$ be given in (\ref{eqe25}), and the sets $\Ub_{\KKI}^{k}$ be defined in (\ref{ep17}) for $k=0, 1,2,3$. Assume, for all $k=0,1,2,3$, we have $\dim(\mathrm{Span}\{\Ub_{\KKI}^{k}\})=3$ and $(h_i^k,r_i^k)$ satisfies the system
	\begin{equation}
	\label{ep26}
	h_i^k a_{v_k}(r_i^k) = \lm_1^k h_i^k + \lm_2^k a_{v_k}(r_i^k)
	\end{equation}
	and
	\begin{equation}
	\label{ep27}
	\frac{(h_i^k)^2}{2}+F_{v_k}(r_i^k) = \lm_1^k r_i^k + \lm_2^k h_i^k
	\end{equation}
	for all $i$ with $\lm_1^k\ne 0$ and $\lm_2^k\ne0$, then $\KKI$ cannot contain a $\Tb_4$.}
\end{a1}

\begin{proof}	
	By Lemma \ref{L100}, it suffices to show that $\ti\Ub_{\KKI}^k:=\lt\{\lt(\begin{matrix}
	h_i^k & r_i^k\\
	a_{v_k}(r_i^k)  & h_i^k\end{matrix}\rt):i=0,1,2,3\rt\}$ cannot \blue{contain} a $\Tb_4$ for some $k\in\{0,1,2,3\}$. Without loss of generality, we assume that $v_0\leq v_1\leq v_2\leq v_3$. \blue{Note that this ordering is only used in the proof of the current lemma.} If $r_i^k=0$ for some $i\ne k$, we have $a_{v_k}(r_i^k)=0$ by (\ref{eq7}). From (\ref{ep26}), it follows that \red{$\lm_1^kh_i^k=0$} and thus \red{$h_i^k=0$}. This means $\ca{\Ub_{\KKI}^{k}}<4$. So we may assume $r_i^k\ne 0$ for all $i\ne k$, and thus
	\begin{equation}\label{ep28}
	v_0<v_1<v_2<v_3.
	\end{equation}
		
	Now we enumerate all possibilities in the following. To simplify notations, we denote
	\begin{equation}\label{ep29}
	D_i^k:=(h_i^k)^2-r_i^k a_{v_k}(r_i^k).
	\end{equation}
	\blue{We may assume that $D_i^k\ne 0$ for all $i\ne k$, as otherwise $\ti\Ub^k_{\KKI}$ would contain Rank-$1$ connections and thus cannot be a $\Tb_4$. This allows us to apply Lemma \ref{l21}.} From (\ref{cgeq21}) it is clear that $h_i^k=-h_k^i$ and $r_i^k=-r_k^i$. By (\ref{eq6}), we calculate
	\begin{eqnarray*}
	a_{v_k}(r_i^k)&\overset{(\ref{eq6})}{=}& \al(v_k+r_i^k)-\al(v_k)\\
	&\overset{(\ref{cgeq21})}{=}&\al(v_i)-\al(v_k)\\
	&=&-\lt(\al(v_k)-\al(v_i)\rt) = -a_{v_i}(r_k^i),
	\end{eqnarray*}
	and thus 
	\begin{equation}\label{ep24}
	D_i^k=D_k^i. 
	\end{equation}
	\red{Note that we have $r_0^3<r_1^3<r_2^3<0$ by (\ref{ep28}), and thus \rrred{by \bblue{(\ref{eq7})}} we have $a_{v_3}(r_0^3)<a_{v_3}(r_1^3)<a_{v_3}(r_2^3)<0$. By Lemma \ref{l20}, we must have $a_{v_3}(r_2^3)>\lm_1^3$. In particular, we must have $\lm_1^3<0$, and thus by Lemma \ref{l21} we have}
\begin{equation}
\label{cgeq26}
	D_2^3<0.
\end{equation}
	
\emph{Case 1.} Assume $\lm_1^0<0$. It follows \red{from (\ref{ep28}) and (\ref{cgeq21})} that $0<r_1^0<r_2^0<r_3^0$ and thus $\lm_1^0<0<a_{v_0}(r_1^0)<a_{v_0}(r_2^0)<a_{v_0}(r_3^0)$. 
By Lemma \ref{l21} and recalling (\ref{ep29}), we have $D_i^0<0$ for $i=1,2,3$. By Proposition \ref{rmk2}, we know that $\ti\Ub^0_{\KKI}$ cannot \blue{contain} a $\Tb_4$.\nl
	
	\emph{Case 2.} Assume $\lm_1^0>0$. By Lemma \ref{l20}, we either have (noting that $h=r=0$ is trivially a solution of (\ref{ep26})-(\ref{ep27}) \rrred{and recalling (\ref{ss63})}) $0<a_{v_0}(r_1^0)<\lm_1^0<a_{v_0}(r_2^0)<a_{v_0}(r_3^0)$ or $0<\lm_1^0<a_{v_0}(r_1^0)<a_{v_0}(r_2^0)<a_{v_0}(r_3^0)$. In the first subcase, by Lemma \ref{l21}, we know $D_1^0>0$. If $D_2^0>0$ and $D_3^0>0$, then $\ti\Ub^0_{\KKI}$ cannot \blue{contain} a $\Tb_4$ by Proposition \ref{rmk2}. So by Lemma \ref{l22}, we only have to consider the cases when $D_2^0>0, D_3^0<0$ or $D_2^0<0, D_3^0<0$. Together with $D_1^0>0$, we are led to two subcases: $D_1^0>0, D_2^0>0, D_3^0<0$ or $D_1^0>0, D_2^0<0, D_3^0<0$. On the other hand, if $0<\lm_1^0<a_{v_0}(r_1^0)<a_{v_0}(r_2^0)<a_{v_0}(r_3^0)$, then by Lemma \ref{l22} and Proposition \ref{rmk2} again, we only have to consider the cases when $D_1^0>0, D_2^0>0, D_3^0<0$ or $D_1^0>0, D_2^0<0, D_3^0<0$. Thus, in conclusion,  we have two subcases to consider. \nl
	
	\emph{Subcase 2.1.} 
	Assume $D_1^0>0, D_2^0>0, D_3^0<0$. By (\ref{ep24}) we have $D_0^3=D_3^0<0$. So if $D_1^3<0$, then by 
\red{(\ref{cgeq26})} we know $D_i^3<0$ for all $i\ne 3$ and thus we are done by Proposition \ref{rmk2}. If $D_1^3>0$, now we have $D_3^1=D_1^3>0$. We claim that $\lm_1^1>0$. Otherwise, we would have $\lm_1^1<0<a_{v_1}(r_3^1)$ and by Lemma \ref{l21} we would have $D_3^1<0$, which is a contradiction. Now as $\lm_1^1>0$, we know from Lemma \ref{l20} that $r_0^1<0$ and $0$ are the only two solutions of the system with $a_{v_1}(r)<\lm_1^1$, and thus $0<\lm_1^1<a_{v_1}(r_2^1)<a_{v_1}(r_3^1)$. Now it follows from Lemma \ref{l22} and $D_3^1>0$ that $D_2^1>0$. As $D_0^1=D_1^0>0$, we have $D_i^1>0$ for $i=0,2,3$, and thus we are done by Proposition \ref{rmk2}. \nl

	\emph{Subcase 2.2.} Assume $D_1^0>0, D_2^0<0, D_3^0<0$. If $D_1^3<0$, we also have $D_0^3=D_3^0<0$ and $D_2^3<0$ by (\ref{cgeq26}). \rrred{So} $D_i^3<0$ for $i=0,1,2$ and we are done by Proposition \ref{rmk2}. If $D_1^3>0$, we have either $D_1^2=D_2^1>0$ or $D_1^2=D_2^1<0$. In the former case, we have $D_i^1>0$ for $i=0,2,3$, and in the latter case \rrred{(recalling (\ref{cgeq26}))} we have $D_i^2<0$ for $i=0,1,3$. Thus in both cases we are done by Proposition \ref{rmk2}. This completes the proof of Lemma \ref{L23}.
\end{proof}

\begin{proof}[Proof of Theorem \blue{\ref{T3}} completed]
	If $\dim\lt(\mathrm{Span}\{\Ub_{\KKI}^{0}\}\rt)=3$, we have two cases: either $\mathrm{Rank}(\AI_{\KKI}^0)=3$ or $\mathrm{Rank}(\AI_{\KKI}^0)=2$ (if $\mathrm{Rank}(\AI_{\KKI}^0)=1$, we must have $\vec h\times \vec r=0$ which is done in Lemma \ref{L1.5}). The latter case is treated in Lemmas \ref{L17} and \ref{L23} (together with explanations immediately before Lemma \ref{L23}), and the former case is treated in Lemmas \ref{L12} and \ref{L15}. 
\end{proof}

\section{Non-existence of two-dimensional $\Tb_4$}

In this section we show that if $\mathrm{dim}\lt(\mathrm{Span}\{\Ub_{\KKI}^{0}\}\rt)=2$ then $\KKI$ cannot \blue{contain} a $\Tb_4$.

%
%

\begin{a2}
	\label{LC3}
	Let $\Ub_{\KKI}^{0}$ be defined by (\ref{eqe25}). If 
	$\mathrm{dim}\lt(\mathrm{Span}\{\Ub_{\KKI}^{0}\}\rt)=2$ then $\KKI$ cannot \blue{contain} a $\Tb_4$.
\end{a2}
\begin{proof}
	By Lemma \ref{L1.22}, we know that $\mathrm{Rank}(\SI_{\KKI}^0)=2$. Using Lemma \ref{L1.5}, we may assume that $\vec h\times\vec r\ne 0$ and $\vec h\times\vec z\ne 0$. \blue{In particular, the first and \rblue{the third} columns in $\SI_{\KKI}^0$ are linearly independent.} So there exist $\gamma_1,\gamma_2$, $\lm_1$, $\lm_2$ and $\mu_1$, $\mu_2$ such that
	\begin{equation}
		\label{eqp4}
		r_i=\gamma_1 h_i + \gamma_2 a(r_i),
	\end{equation}
	\begin{equation}
		\label{eqp5}
		h_i a(r_i)=\lm_1 h_i + \lm_2 a(r_i),
	\end{equation}
	and
	\begin{equation}
		\label{eqp6}
		\frac{h_i^2}{2}+F(r_i)=\mu_1 h_i + \mu_2 a(r_i).
	\end{equation}
	It follows that
	\begin{equation*}
		\mathrm{Span}\{\Ub_{\KKI}^{0}\}\overset{(\ref{eq400}),(\ref{eqe25})}{=}\lt\{\lt(\begin{array}{cc} s & \gamma_1 s+\gamma_2 t \\ t & s \\ \lm_1 s+\lm_2 t & \mu_1 s + \mu_2 t \end{array}\rt): s,t\in\R\rt\}.
	\end{equation*}
	The three minors in $\mathrm{Span}\{\Ub_{\KKI}^{0}\}$ are
	\begin{equation}
		\label{eq529aabb11}
		M_1 = s^2 - \gamma_1 st - \gamma_2 t^2,
	\end{equation}
	\begin{equation*}
		\begin{split}
			M_2 &= s\lt(\mu_1 s + \mu_2 t\rt) - \lt(\gamma_1 s+\gamma_2 t\rt)\lt(\lm_1 s+\lm_2 t\rt)\\
			&= \lt(\mu_1-\gamma_1\lm_1\rt)s^2 + \lt(\mu_2-\gamma_1\lm_2-\gamma_2\lm_1\rt)st -\gamma_2\lm_2 t^2,
		\end{split}
	\end{equation*}
	and
	\begin{equation*}
		\begin{split}
			M_3 &= t\lt(\mu_1 s + \mu_2 t\rt) - s\lt(\lm_1 s+\lm_2 t\rt)\\
			&= -\lm_1 s^2 + \lt(\mu_1-\lm_2\rt)st + \mu_2 t^2. 
		\end{split}
	\end{equation*}
	
	If $\gamma_1^2+4\gamma_2<0$, then \rrred{(viewing the left hand side as a quadratic \rblue{polynomial} in $s$)}
	\begin{equation*}
		s^2 - \gamma_1 st -\gamma_2 t^2>0 
	\end{equation*}
	for all $(s,t)\ne (0,0)$ and \rred{ so we see from (\ref{eq529aabb11})} that $\mathrm{Span}\{\Ub_{\KKI}^{0}\}$ has no Rank-$1$ directions. If $\gamma_1^2+4\gamma_2=0$, \green{ then 
		$M_1=\lt(s-\frac{\gamma_1 t}{2}\rt)^2$.
		So $s=\frac{\gamma_1 t}{2}$} produces the only possible Rank-$1$ direction in $\mathrm{Span}\{\Ub_{\KKI}^{0}\}$ and we can apply Lemma \ref{L13} (a). So for the rest of the proof we assume that $\gamma_1^2+4\gamma_2>0$, which implies that the equation \blue{$\blue{x}^2 - \gamma_1 \blue{x} -\gamma_2=0$ has two distinct solutions and thus one can write $\blue{x}^2 - \gamma_1 \blue{x} -\gamma_2=(x-k)(x-l)$ for some $k\ne l$. It follows that $\frac{s^2}{t^2}-\gamma_1\frac{s}{t}-\gamma_2=(\frac{s}{t}-k)(\frac{s}{t}-l)$ and therefore}
	\begin{equation}
		\label{eqp8}
		s^2 - \gamma_1 st -\gamma_2 t^2=\lt(s-kt\rt)\lt(s-lt\rt).
	\end{equation}
	The Rank-$1$ directions in $\mathrm{Span}\{\Ub_{\KKI}^{0}\}$ require $M_1 = M_2 = M_3=0$. From (\ref{eqp8}), the only possible Rank-$1$ directions in $\mathrm{Span}\{\Ub_{\KKI}^{0}\}$ must satisfy $s=kt$ or $s=lt$. Now we check these two directions.
	
	Note that from (\ref{eqp8}), we have 
	\begin{equation}
		\label{eqp9}
		\gamma_1=k+l, \qd\gamma_2=-kl.
	\end{equation} 
	When $s=kt$, plugging this into $M_2$ and $M_3$ and using (\ref{eqp9}) give
	\begin{equation}
		\label{eqp10}
		\begin{split}
			M_2 &= \lt(\mu_1-\gamma_1\lm_1\rt)k^2t^2 + \lt(\mu_2-\gamma_1\lm_2-\gamma_2\lm_1\rt)kt^2 -\gamma_2\lm_2 t^2\\
			&\overset{\rrred{(\ref{eqp9})}}{=} \lt(\mu_1k^2-(k+l)\lm_1k^2+\mu_2 k-(k+l)\lm_2k+kl\lm_1k+kl\lm_2\rt)t^2\\
			&= \lt(-\lm_1 k^2+\lt(\mu_1-\lm_2\rt)k+\mu_2\rt)kt^2
		\end{split}
	\end{equation}
	and
	\begin{equation}
		\label{eqp11}
		\begin{split}
			M_3 &= -\lm_1 k^2 t^2 + \lt(\mu_1-\lm_2\rt)kt^2 + \mu_2 t^2\\
			&= \lt(-\lm_1 k^2+(\mu_1-\lm_2)k+\mu_2\rt)t^2.
		\end{split}
	\end{equation}
	When $s=lt$, \blue{with $k$ and $l$ switched in (\ref{eqp10}) and (\ref{eqp11}) we obtain}
	\begin{equation*}
		M_2 = \lt(-\lm_1 l^2+\lt(\mu_1-\lm_2\rt)l+\mu_2\rt)lt^2
	\end{equation*}
	and
	\begin{equation*}
		M_3 = \lt(-\lm_1 l^2+(\mu_1-\lm_2)l+\mu_2\rt)t^2.
	\end{equation*}
	So in order for $s=kt$ and $s=lt$ to be Rank-$1$ directions in $\mathrm{Span}\{\Ub_{\KKI}^{0}\}$, we need $k\ne l$ to satisfy the quadratic equation
	$-\lm_1 x^2+(\mu_1-\lm_2)x+\mu_2=0$. As $k\ne l$ also solve $x^2 - \gamma_1 x -\gamma_2=0$, it follows that
	\begin{align*}
		\mu_1-\lm_2&=\lm_1\gamma_1,\nn\\
		\mu_2&=\lm_1\gamma_2.
	\end{align*}
	Now rewriting \eqref{eqp6} with $\mu_1=\lm_1\gamma_1+\lm_2$ and $\mu_2=\lm_1\gamma_2$, we obtain
	\begin{equation*}
		\frac{h_i^2}{2}+F(r_i)=\lt(\lm_1\gamma_1+\lm_2\rt) h_i + \lm_1\gamma_2 a(r_i).
	\end{equation*}
	From \eqref{eqp4}, we write $\gamma_2 a(r_i)=r_i-\gamma_1h_i$. Plugging this into the above gives
	\begin{align*}
		\frac{h_i^2}{2}+F(r_i)&=\lt(\lm_1\gamma_1+\lm_2\rt) h_i + \lm_1\lt(r_i-\gamma_1h_i\rt)\\
		&=\lm_1 r_i+\lm_2 h_i.
	\end{align*}
This together with \eqref{eqp5} shows that $(h_i, r_i)$ satisfies
\begin{align*}
	h_i a(r_i)&=\lm_1 h_i + \lm_2 a(r_i),\\
	\frac{h_i^2}{2}+F(r_i)&=\lm_1 r_i+\lm_2 h_i,
\end{align*}
which is exactly the system \eqref{eqp26}-\eqref{eqp27}. 

Now recall $h_i^k, r_i^k$ and the sets $\Ub_{\KKI}^{k}$ defined in \eqref{cgeq21} and \eqref{ep17}, respectively. From the proof of Lemma \ref{L1.2} it is clear that $\mathrm{dim}\lt(\mathrm{Span}\{\Ub_{\KKI}^{k}\}\rt)=2$ for all $k=0, 1, 2, 3$. Applying Lemma \ref{L1.2} to $\Ub_{\KKI}^{k}$, it suffices to show that $\Ub_{\KKI}^{k}$ cannot contain a $\Tb_4$ for some $k$. \nl

\emph{Case 1.} If for some $k\in\{0, 1, 2, 3\}$, there does not exist $\lambda_1^k, \lambda_2^k$ such that \eqref{ep26}-\eqref{ep27} are satisfied, then the above calculations show that $\mathrm{Span}\{\Ub_{\KKI}^{k}\}$ contains at most one Rank-$1$ direction, and thus $\Ub_{\KKI}^{k}$ cannot contain a $\Tb_4$. \nl

\emph{Case 2.} Assume for all  $k=0, 1, 2, 3$, there exist $\lambda_1^k, \lambda_2^k$ such that \eqref{ep26}-\eqref{ep27} are satisfied. Recall the sets $\ti\Ub_{\KKI}^k:=\lt\{\lt(\begin{matrix}
	h_i^k & r_i^k\\
	a_{v_k}(r_i^k)  & h_i^k\end{matrix}\rt):i=0,1,2,3\rt\}$. The linear mapping $L$ defined in the proof of Lemma \ref{L100} with index $0$ replaced by $k$ still satisfies \eqref{eqa30}, and is an isomorphism because $\mathrm{dim}\lt(\mathrm{Span}\{\ti\Ub_{\KKI}^k\}\rt)=2$ (as a result of e.g. $\vec{h^k}\times \vec{r^k} \ne 0$). By Lemma \ref{LL0}, it suffices to show that $\ti\Ub_{\KKI}^k$ cannot contain a $\Tb_4$ for some $k$. Note that in the proofs of Lemmas \ref{L17} and \ref{L23}, the assumptions $\mathrm{dim}\lt(\mathrm{Span}\{\Ub_{\KKI}^0\}\rt)=3$ and $\mathrm{dim}\lt(\mathrm{Span}\{\Ub_{\KKI}^k\}\rt)=3$ are only used to ensure that we can apply Lemma \ref{L100}. In particular, the essential parts of the proofs, i.e. $\ti\Ub_{\KKI}^0$ not containing $\Tb_4$ in Lemma \ref{L17} and $\ti\Ub_{\KKI}^k$ not containing $\Tb_4$ in Lemma \ref{L23}, do not use the dimension $3$ assumption. If for some $k$, we have $\lambda_1^k=0$ or $\lambda_2^k=0$, then the proof of Lemma \ref{L17} shows that $\ti\Ub_{\KKI}^k$ cannot contain a $\Tb_4$. Finally, if $\lambda_1^k\ne 0$ and $\lambda_2^k\ne0$ for all $k$, then the proof of Lemma \ref{L23} shows that $\ti\Ub_{\KKI}^k$ cannot contain a $\Tb_4$ for some $k\in\{0, 1, 2, 3\}$. This concludes the proof.
\end{proof}

Finally, putting Theorems \ref{T3} and \ref{LC3} together, we complete the proof of Theorem \ref{T1}.

\section{Proof of Proposition \ref{TCXE}}\label{ce}

\zred{We start by giving} a more explicit equivalent condition for the set $\K_1$ to contain Rank-$1$ connections.

\begin{a1}\label{L25}
\zred{Let $I$ be an interval, and let the set $\K^I_1$ be defined in (\ref{eq4039}) with the function $\al\in C^2(\R)$ satisfying $\al'>0$. Then the set $\K^I_1$ contains Rank-$1$ connections if and only if there exist $v\in\bbblue{I}$ and $r\ne 0$ such that \bbblue{$v+r\in I$ and}
	\begin{equation}\label{ep105}
	2F_v(r)=ra_v(r),
	\end{equation}
	where the functions $a_v$ and $F_v$ are defined in (\ref{eq6}).}
\end{a1}

\begin{proof}
	By definition, the set $\zred{\K^I_1}$ contains Rank-$1$ connections if and only if there exist $(u,v)\ne(\ti u,\ti v)$ such that $\zred{v, \ti{v}\in I}$ and $\mathrm{Rank}\lt(P(\ti u,\ti v)-P(u,v)\rt)=1$, where the mapping $P$ is given in (\ref{eq4}). Denoting by $h=\ti u-u$, $r=\ti v-v$ and recalling the notations in (\ref{eq6}) and (\ref{eq400}), it follows from Lemma \ref{L5} that there exists an invertible matrix $B$ such that $B(P(\ti u,\ti v)-P(u,v))=\QI_v(h,r)$, where $\QI_v(h,r)$ is given in (\ref{eq400}). Hence $\mathrm{Rank}\lt(P(\ti u,\ti v)-P(u,v)\rt)=1$ if and only if $\mathrm{Rank}\lt(\QI_v(h,r)\rt)=1$. Therefore the set $\zred{\K^I_1}$ contains Rank-$1$ connections if and only if there exist $v\in I$ and $(h,r)\ne (0,0)\in\R^2$ such that $v+r\in I$ and $\mathrm{Rank}\lt(\QI_v(h,r)\rt)=1$.
	
	Given $v\in\R$ and $(h,r)\ne (0,0)$, we claim that $\mathrm{Rank}\lt(\QI_v(h,r)\rt)=1$ if and only if
	\begin{equation}\label{ep106}
	h^2 = ra_v(r)\qd\text{ and }\qd 2F_v(r)=ra_v(r).
	\end{equation}
	To see this, we write out the three minors of $\QI_v(h,r)$:
	\begin{equation*}
	M_1 = h^2 - ra_v(r), \qd M_2=\frac{h^3}{2}+hF_v(r)-rha_v(r),
	\end{equation*}
	and
	\begin{equation*}
	M_3=\frac{h^2}{2}a_v(r)+a_v(r)F_v(r)-h^2a_v(r).
	\end{equation*}
	If $\mathrm{Rank}\lt(\QI_v(h,r)\rt)=1$, then $M_1=M_2=M_3=0$. From $M_1=0$ we obtain $h^2 = ra_v(r)$. Note that from this, (\ref{eq7}) and $(h,r)\ne(0,0)$, we must have $h\ne 0$, $r\ne 0$ and $a_v(r)\ne 0$. Now $M_2=0$ and $M_3=0$ reduce to 
	\begin{equation}\label{ep107}
	\frac{h^2}{2}+F_v(r)-ra_v(r)=0 \qd\text{ and }\qd \frac{h^2}{2}+F_v(r)-h^2=0.
	\end{equation}
	Comparing the equations in (\ref{ep107}) and substituting $h^2$ by $ra_v(r)$, one readily sees that  (\ref{ep107}) is equivalent to (\ref{ep106}). Conversely, if (\ref{ep106}) holds, then we have (\ref{ep107}) and it is clear that $M_1=M_2=M_3=0$. Thus we have $\mathrm{Rank}\lt(\QI_v(h,r)\rt)=1$.
	
	Now if $\zred{\K^I_1}$ contains Rank-$1$ connections, then there exist $v\in I$ and $(h,r)\ne (0,0)\in\R^2$ such that $v+r\in I$ and $\mathrm{Rank}\lt(\QI_v(h,r)\rt)=1$. Therefore (\ref{ep106}) and thus (\ref{ep105}) hold true. Conversely, if (\ref{ep105}) holds for some $v$ and $r\ne 0$, then as $a_{v}(0)=0$ and $a_v'>0$ (recalling (\ref{eq7})), it is clear that $ra_v(r)>0$ and thus one can choose $h=\sqrt{ra_v(r)}$. With this choice of $v, r, h$, the equations in (\ref{ep106}) are satisfied. Hence $\mathrm{Rank}\lt(\QI_v(h,r)\rt)=1$ and $\zred{\K^I_1}$ contains Rank-$1$ connections.
\end{proof}


\begin{proof}[Proof of Proposition \ref{TCXE}]
\zred{First we assume} that $\al$ has an isolated inflection point at \zred{$v_0\in I$}. Without loss of generality, assume that 
\begin{equation}\label{ep112}
\al''(v)<0\qd\text{ for } v_0-\delta<v<v_0
\end{equation}
and 
\begin{equation}\label{ep113}
\al''(v)>0\qd\text{ for } v_0<v<v_0+\delta
\end{equation}
for some $\delta>0$ \bbblue{sufficiently small}. Recall the definitions of the translation functions $a_v$ and $F_v$ in (\ref{eq6}) and the properties listed in (\ref{eq7})-(\ref{eq8}). As in the proof of Lemma \ref{l18}, we define
\begin{equation}
\label{ep3010}
g_v(r):= 2F_v(r)-ra_v(r)
\end{equation}
and obtain
\begin{equation}\label{ep115}
g_v(0)=g_v'(0)=0,\qd g_v''(r)=-ra_v''(r).
\end{equation}
By Lemma \ref{l18}, since $F_v'=a_v$ and $a_v'(t)\overset{(\ref{eq6})}{=}\al'(v+t)$ we have 
\begin{equation*}
\begin{split}
&\text{if for some }\ggreen{r_0, r_1>0}\text{ we have }\al''>0\text{ in }\ggreen{\lt(v-r_0,v+r_1\rt)}, \text{ then }\\
&\qd\qd g_v(r)>0\text{ for }r\in (-r_0,0), \qd g_v(r)<0\text{ for }r\in (0,\ggreen{r_1}),
\end{split}
\end{equation*}
and
\begin{equation}\label{ep117}
\begin{split}
&\text{if for some }\ggreen{r_0, r_1>0}\text{ we have }\al''<0\text{ in }\ggreen{\lt(v-r_0,v+r_1\rt)}, \text{ then }\\
&\qd\qd g_v(r)<0\text{ for }r\in (-r_0,0), \qd g_v(r)>0\text{ for }r\in (0,\ggreen{r_1}).
\end{split}
\end{equation}
	
Next we define the functions
\begin{equation*}
p(v):=g_v(v_0-v)=2F_v(v_0-v)-(v_0-v)a_v(v_0-v)
\end{equation*}
and
\begin{equation*}
q(v):=g_v\lt(v_0+\frac{\delta}{2}-v\rt)=2F_v\lt(v_0+\frac{\delta}{2}-v\rt)-\lt(v_0+\frac{\delta}{2}-v\rt)a_v\lt(v_0+\frac{\delta}{2}-v\rt).
\end{equation*}
Using the definitions for $a_v$ and $F_v$ as in (\ref{eq6}), we write out 
\begin{equation*}
p(v)=2\lt(\Fk(v_0)-\Fk(v)-\al(v)\lt(v_0-v\rt)\rt)-(v_0-v)\lt(\al(v_0)-\al(v)\rt)
\end{equation*}
and clearly $p(v)$ is continuous. Similarly $q(v)$ is also continuous. It follows from (\ref{ep112}) and (\ref{ep117}) that 	
\begin{equation}\label{ep121}
p(v) = g_v(v_0-v) >0 \qd\text{ for all } v\in \lt(v_0-\frac{\delta}{2},v_0\rt).
\end{equation}
On the other hand, note that $q(v_0)=g_{v_0}(\frac{\delta}{2})$. We deduce from (\ref{ep115}) and (\ref{ep112})-(\ref{ep113}) that $g_{v_0}(0)=g_{v_0}'(0)=0$ and $g_{v_0}$ is concave locally around the origin, and thus $q(v_0)=g_{v_0}(\frac{\delta}{2})<0$. As $q$ is continuous, it follows immediately that
\begin{equation}\label{ep122}
q(v_1)=g_{v_1}\lt(v_0+\frac{\delta}{2}-v_1\rt)<0 \qd\text{ for some } v_1\in \lt(v_0-\frac{\delta}{2},v_0\rt). 
\end{equation}
\zred{Consider the continuous function
\begin{equation}
\label{eq3004}
\varpi(w)=g_{v_1}(v_0+w-v_1)\qd\text{ for }w\in \lt[0,\frac{\delta}{2}\rt].
\end{equation}
Note that since $v_1\in \lt(v_0-\frac{\delta}{2},v_0\rt)$ we have that $\varpi(0)\overset{(\ref{eq3004}), (\ref{ep121})}{=}p(v_1)>0$ and $\varpi\lt(\frac{\delta}{2}\rt)\overset{(\ref{eq3004}), (\ref{ep122})}{=}q(v_1)<0$. So 
there exists $v_2\in(v_0,v_0+\frac{\delta}{2})$ such that $\varpi(v_2-v_0)\overset{(\ref{eq3004})}{=}g_{v_1}(v_2-v_1)=0$}. Denoting by $r:=v_2-v_1>0$, this translates to $2F_{v_1}(r)-ra_{v_1}(r)\overset{(\ref{ep3010})}{=}0$, and thus gives Rank-$1$ connection in the set $\K_1^I$ by Lemma \ref{L25}.
	
\zred{Now suppose $\al$ is either strictly convex or strictly concave on $I$, then by Lemma \ref{l18} and Lemma \ref{L25}, the set $\mathcal{K}_1^I$ contains no Rank-$1$ connections.}
\end{proof}


\begin{thebibliography}{99}

	
	
	

	

	
	\bibitem[Bi-Br 05]{bres}   S. Bianchini; A. Bressan. \emph{Vanishing viscosity solutions of nonlinear hyperbolic systems}. Ann. of Math. (2) 161 (2005), no. 1, 223--342.

\bibitem[Br-Cr-Pi 00]{bres2}   A. Bressan; G. Crasta; B. Piccoli. 
\emph{Well-posedness of the Cauchy problem for} $n\times n$ \emph{ systems of conservation laws}.  Mem. Amer. Math. Soc. 146 (2000), no. 694.
	\bibitem[Bu-De-Is-Sz 15]{buvi2} 
T. Buckmaster; C. De Lellis; P. Isett; L. Sz\'{e}kelyhidi, Jr. \emph{Anomalous dissipation for $1/5$-H\"{o}lder Euler flows.} Ann. of Math. (2) 182 (2015), no. 1, 127--172. 



\bibitem[Bu-De-Sz-Vi 19]{camles8}   
T. Buckmaster; C. De Lellis; L. Sz\'{e}kelyhidi, Jr.; V. Vicol. \emph{Onsager's conjecture for admissible weak solutions.} Comm. Pure Appl. Math. 72 (2019), no. 2, 229--274.

\bibitem[Bu-Vi 19]{buvi}   
\blue{T. Buckmaster; V. Vicol. \emph{Nonuniqueness of weak solutions to the Navier-Stokes equation.} Ann. of Math. (2) 189 (2019), no. 1, 101--144.}
	
\bibitem[Ch-Ki 02]{chki}  M. Chlebik; B. Kirchheim.  \emph{Rigidity for the four gradient problem.}  J. Reine Angew. Math. 551 (2002), 1--9.  

\bibitem[De-De-Ki-Ti 19]{dede}  C. De Lellis; G. De Philippis; B. Kirchheim; R. Tione. \xred{\emph{Geometric measure theory and differential inclusions}.   Preprint. https://arxiv.org/abs/1910.00335}


\bibitem[De-Sz 09]{camles1} 
	C. De Lellis; L. Sz\'{e}kelyhidi. \emph{The Euler equations as a differential inclusion.} Ann. of Math. (2) 170 (2009), no. 3, 1417--1436.

\bibitem[De-Sz 12]{camles3} 
	C. De Lellis; L. Sz\'{e}kelyhidi. \emph{The h-principle and the equations of fluid dynamics}. Bull. Amer. Math. Soc. (N.S.) 49 (2012), no. 3, 347--375.
 
\bibitem[De-Sz 13]{camles2} 
	C. De Lellis; L. Sz\'{e}kelyhidi. \emph{Dissipative continuous Euler flows}. Invent. Math. 193 (2013), no. 2, 377--407.



\bibitem[De-Sz 19]{camlas10} 
	C. De Lellis; L. Sz\'{e}kelyhidi. \emph{On turbulence and geometry: from Nash to Onsager}. To appear in the Notices of the AMS.
	
	
	
	
	

	\bibitem[DP 83]{dp2}
	R. J. DiPerna. \emph{Convergence of approximate solutions to conservation laws.} Arch. Rational Mech. Anal. 82 (1983), no. 1, 27--70.
	
	\bibitem[DP 85]{dp}
	R. J. DiPerna. \emph{Compensated compactness and general systems of conservation laws.} Trans. Amer. Math. Soc. 292 (1985), no. 2, 383--420.
	
	
	
	\bibitem[Ev 86]{evans} 
	L. C. Evans. \emph{Quasiconvexity and partial regularity in the calculus of variations}. Arch. Rational Mech. Anal. 95 (1986), no. 3, 227--252. 

	\bibitem[Ev 90]{evans2} 
	L. C. Evans. \emph{Weak Convergence Methods in Partial Differential Equations}. CBMS Regional Conference Series in Mathematics, 74. Published for the Conference Board of the Mathematical Sciences, Washington, DC; by the American Mathematical Society, Providence, RI, 1990. 

	\bibitem[Ev 10]{evans3} L. C. Evans. \emph{Partial differential equations.} Second edition. Graduate Studies in Mathematics, 19. American Mathematical Society, Providence, RI, 2010.

	\bibitem[Fa-Sz 08]{faclas} D. Faraco; L. Sz\'{e}kelyhidi. \emph{Tartar's conjecture and localization of the quasiconvex hull in $\R^{2\times 2}$}.  Acta Math. 200 (2008), no. 2, 279--305.

\bibitem[F\"{o}-Sz 18]{fosz}  C. F\"{o}rster;  L. Sz\'{e}kelyhidi. \emph{$T_5$-configurations and non-rigid sets of matrices}. Calc. Var. Partial Differential Equations 57 (2018), no. 1, Art. 19, 12 pp.
	\bibitem[Gr 86]{grom} 
	M. Gromov. \emph{Partial differential relations.} Ergebnisse der Mathematik und ihrer Grenzgebiete (3) [Results in Mathematics and Related Areas (3)], 9. Springer-Verlag, Berlin, 1986.
	
		\bibitem[Is 13]{phil2} P. Isett. \emph{H\"{o}lder continuous Euler flows with compact support in time}. Thesis (Ph.D.)--Princeton University. 2013. 227 pp. 

	\bibitem[Is 17]{phil}   
	P. Isett. \emph{H\"{o}lder continuous Euler flows in three dimensions with compact support in time.} Annals of Mathematics Studies, 196. Princeton University Press, Princeton, NJ, 2017. 
	 
	\bibitem[Is 18]{phil1} P. Isett. \emph{A proof of Onsager's conjecture}. Ann. of Math. (2) 188 (2018), no. 3, 871--963. 


	   

	

	
	
	
	\bibitem[Ki 03]{kir}
	B. Kirchheim. \emph{Rigidity and geometry of microstructures.} Habilitation Thesis, University of Leipzig, 2003. MIS.MPG preprint 16/2003.
	
	\bibitem[Ki-M\"{u}-\v{S}v 03]{kms}
	B. Kirchheim; S. M\"{u}ller; V. \v{S}ver\'{a}k. \emph{Studying nonlinear pde by geometry in matrix space.} Geometric analysis and nonlinear partial differential equations, 347--395, Springer, Berlin, 2003.
	
	
	\bibitem[Ku 55]{kuip} N. H. Kuiper. \emph{On $C^1$-isometric imbeddings.} I, II. Nederl. Akad. Wetensch. Proc. Ser. A. 58 = Indag. Math. 17 (1955), 545--556, 683--689.

	\bibitem[Lo-Pe 19]{lope} 
	A. Lorent; G. Peng.  \emph{Null Lagrangian Measures in subspaces, compensated compactness and conservation laws.}  Arch. Ration. Mech. Anal. 234 (2019), no. 2, 857--910.
	

	
	\bibitem[M\"{u} 99]{mul}
	S. M\"{u}ller. \emph{Variational models for microstructure and phase transitions.} Calculus of variations and geometric evolution problems (Cetraro, 1996), 85--210, Lecture Notes in Math., 1713, Fond. CIME/CIME Found. Subser., Springer, Berlin, 1999.
	
	\bibitem[M\"{u}-Ri-\v{S}v 05]{mulsvri} 
	S. M\"{u}ller; M. O. Rieger; \blue{V. \v{S}ver\'{a}k}. \emph{Parabolic systems with nowhere smooth solutions}.  Arch. Ration. Mech. Anal. 177 (2005), no. 1, 1--20.
\bibitem[M\"{u}-\v{S}v 96]{mulsv1} S. M\"{u}ller; V. \v{S}ver\'{a}k. \emph{Attainment results for the two-well problem by convex integration.} Geometric analysis and the calculus of variations, 239--251, Int. Press, Cambridge, MA, 1996.

\bibitem[M\"{u}-\v{S}v 99]{mulsv3} 
S. M\"{u}ller; V. \v{S}ver\'{a}k.  \emph{Convex integration with constraints and applications to phase transitions and partial differential equations}. J. Eur. Math. Soc. (JEMS) 1 (1999), no. 4, 393--422.

\bibitem[M\"{u}-\v{S}v 03]{mulsv2} 
S. M\"{u}ller; V. \v{S}ver\'{a}k.  \emph{Convex integration for Lipschitz mappings and counterexamples to regularity.} Ann. of Math. (2) 157 (2003), no. 3, 715--742.


	


\bibitem[Mu 78]{murat}  F. Murat. \emph{Compacite par compensation}. Ann. Scuola Norm. Sup. Pisa Sci. Fis.
Mat. 5, 489--507 (1978).

\bibitem[Na 54]{nash} 
J. Nash.  \emph{$C^1$  isometric imbeddings.} Ann. of Math. (2) 60, (1954). 383-396. 
	
\bibitem[Ol 57]{ol}  O. A. Oleinik. \emph{Discontinuous solutions of nonlinear differential equations}. Usp. Mat. Nauk.
12 (1957), pp. 3--73; English transl. in AMS Transl. 26 (1963), pp. 1155--1163. 

  
\bibitem[Sc 74]{sch} V. Scheffer. \emph{Regularity and irregularity of solutions to nonlinear second order
elliptic systems of partial differential equations and inequalities.} Thesis (Ph.D.)--Princeton University. 1974. 116 pp. 

\bibitem[Sc 93]{sch2} V. Scheffer. \emph{An inviscid flow with compact support in space-time.} J. Geom. Anal. 3 (1993), no. 4, 343–401.

\bibitem[Sh 97]{shnir} A. Shnirelman. \emph{ On the nonuniqueness of weak solution of the Euler equation.} Comm. Pure Appl. Math. 50 (1997), no. 12, 1261--1286.

	
	

\bibitem[\v{S}v 92]{sv5} 
V. \v{S}ver\'{a}k.  \emph{New examples of quasiconvex functions.}   Arch. Rational Mech. Anal. 119 (1992), no. 4, 293--300.
	
\bibitem[\v{S}v 16]{svper} V. \v{S}ver\'{a}k.  \emph{Personal communication.} Sabbatical visit. Minnesota, 2016.

\bibitem[\v{S}v 18]{svper2} V. \v{S}ver\'{a}k.  \emph{Personal communication.} Research visit. Minnesota, 2018.

	


\bibitem[Sz 04]{sz2} L. Sz\'{e}kelyhidi, Jr. \emph{The regularity of critical points of polyconvex functionals.} Arch. Ration. Mech. Anal. 172 (2004), no. 1, 133--152.
\bibitem[Sz 05]{sz} L. Sz\'{e}kelyhidi, Jr. \emph{ Rank-one convex hulls in $\R^{2\times2}$.} Calc. Var. Partial Diff. Eq. 22 (2005)	no. 3, 253--281.	

\bibitem[Ta 79]{ta1} L. Tartar. \emph{Compensated compactness and applications to partial differential equations.} Nonlinear analysis and mechanics: Heriot-Watt Symposium, Vol. IV, pp. 136--212, Res. Notes in Math., 39, Pitman, Boston, Mass.-London, 1979.
	
\bibitem[Ta 83]{ta2} L. Tartar. \emph{The compensated compactness method applied to systems of conservation laws.} Systems of nonlinear partial differential equations (Oxford, 1982), 263--285, NATO Adv. Sci. Inst. Ser. C Math. Phys. Sci., 111, Reidel, Dordrecht, 1983.

\bibitem[Ta 93]{ta3} \rblue{L. Tartar. \emph{Some remarks on separately convex functions.} Microstructure and phase transition, 191--204, IMA Vol. Math. Appl., 54, Springer, New York, 1993.}
	
	
\end{thebibliography}
\end{document}